\documentclass[times,final]{elsarticle}


\usepackage{framed,multirow}

\usepackage{amssymb}
\usepackage{latexsym}
\usepackage{mathrsfs}
\usepackage{url}
\usepackage{xcolor}
\definecolor{newcolor}{rgb}{.8,.349,.1}

\RequirePackage{xcolor}
\usepackage{algorithm}
\usepackage{algorithmic}
\usepackage{amsfonts}
\usepackage{amsmath}
\usepackage{array}
\usepackage[T1]{fontenc}
\usepackage{hyperref}
\hypersetup{colorlinks,linkcolor=black,urlcolor=black,citecolor=black}
\usepackage{multirow}
\usepackage{tikz}
\usepackage{bm}
\newtheorem{remark}{Remark}

\usepackage{subfigure}

\usepackage{graphicx}
\usepackage{epstopdf}

\begin{document}

\begin{frontmatter}

\title{Subspace method based on neural networks for solving eigenvalue problems \tnoteref{tnote1}}%
\tnotetext[tnote1]{This work was supported by  the National Natural Science Foundation of China under grants 92270206 and 12071045 and the Fund of National Key Laboratory of Computational Physics }

\author[1]{Xiaoying Dai\corref{cor1}}
\cortext[cor1]{Corresponding author:
  email: daixy@lsec.cc.ac.cn;}
\author[1]{Yunying Fan}
\ead{fanyunying@amss.ac.cn}
\author[2]{Zhiqiang Sheng}
\ead{sheng\_zhiqiang@iapcm.ac.cn}

\address[1]{SKLMS, Academy of Mathematics and Systems Science, Chinese Academy of Sciences, Beijing 100190, China; and School of Mathematical Sciences, University of Chinese Academy of Sciences, Beijing 100049, China}
\address[2]{National Key Laboratory of Computational Physics, Institute of Applied Physics and Computational Mathematics, Beijing, 100088, China, and HEDPS, Center for Applied Physics and Technology, and College of Engineering, Peking University, Beijing, 100871, China}

\begin{abstract}
In this paper, we propose a subspace method based on neural networks for eigenvalue problems with high accuracy and low cost. We first construct a neural network-based orthogonal basis by some deep learning method and dimensionality reduction technique, and then calculate the Galerkin projection of the eigenvalue problem onto the subspace spanned by the orthogonal basis and obtain an approximate solution. Numerical experiments show that we can obtain approximate eigenvalues and eigenfunctions with very high accuracy but low cost. 
\end{abstract}


\end{frontmatter}



\section{Introduction}
\label{intro}

Eigenvalue problems play an important role in various fields of mathematics, physics, and chemistry, such as material mechanics, fluid mechanics, structural mechanics, and quantum mechanics. Accordingly, the study of efficient numerical discretization methods to tackle these problems is an active field in scientific computing. Over recent decades, many numerical methods have been developed to discretize eigenvalue problems with great success, including the finite difference method \cite{kuttler1984eigenvalues}, the finite volume method \cite{dai2011finite,liang2001finite}, the finite element method \cite{BABUSKA1991641,boffi2010finite,sun2016finite}, and the spectral method \cite{feit1982solution}. 

With the rapid development of machine learning, many numerical discretization methods for partial differential equations (PDEs) based on deep neural networks have received increasing attention, such as the physical-information neural networks (PINN) \cite{raissi2019physics}, the deep Galerkin method (DGM) \cite{sirignano2018dgm}, the deep Ritz method (DRM) \cite{yu2018deep}, and the weak adversarial networks (WAN) \cite{2020wan}. By these methods, the solution of a partial differential equation is represented by {neural networks directly} and the optimal parameters of neural networks are obtained by minimizing some loss function. The main difference between these methods lies in the choice of the loss function. For instance, PINN and DGM use the $L^2$ norm of the residual in the strong form as the loss function, while DRM and WAN use the energy functional and the operator norm in the weak form as the loss function, respectively. 

Due to the powerful approximation capabilities of neural networks \cite{2024nn,devore2021neural,he2018relu,he2023deep,LESHNO1993861,lu2021deep,ma2022barron,shen2019deep,siegel2020approximation,siegel2024sharp}, the aforementioned methods are used to deal with many different kinds of partial differential equations, including  not only some  source problems but also some eigenvalue problems. For instance, in \cite{Eweinna2019}, a wave function ansatz based on deep neural networks is introduced to solve the many-electron Schr\"{o}dinger equation. In \cite{Hermann_Schätzle_Noé_2020,Hermann_Spencer_Choo_Mezzacapo_Foulkes_Pfau_Carleo_Noé_2023}, various wave function ansatzes in the form of a Slater determinant are also introduced for the Schr\"{o}dinger equation. In \cite{PhysRevA.}, a neural network-based multistate solver is proposed for finding the first $k$ eigenpairs of the Schr\"{o}dinger equation simultaneously. In \cite{Han_Lu_Zhou_2020}, a diffusion Monte Carlo-like method is proposed for handling high-dimensional eigenvalue problems. In \cite{Li_Ying_2022}, a semigroup method is provided for solving high-dimensional elliptic partial differential equations using neural networks and is extended to a primal-dual algorithm for solving eigenvalue problems. In \cite{li2024tensor}, a type of tensor neural network-based on the DECOMP/PARAFAC (CP) tensor decomposition is proposed, along with an efficient numerical integration scheme for tensor-type functions. Theoretical analyses and applications of tensor neural networks are provided in \cite{wang2024solving} and \cite{wang2022solving,WANG2024112928}, respectively. For example, tensor neural networks have been used to solve the smallest eigenvalue and the corresponding eigenfunction of Schr\"{o}dinger equation in \cite{wang2022solving} and the first $k$ eigenpairs of high-dimensional eigenvalue problems in \cite{WANG2024112928}. Additionally, in \cite{yu2024generalizationerrorestimatesmachine}, a machine learning method is proposed for computing eigenvalues and eigenfunctions of the Schr\"{o}dinger equation, and an explicit convergence rate of the generalization error is also provided in the article.

While deep neural network-based methods have made strides in solving high dimensional problems where classical methods fail, they face challenges in solving low-dimensional problems due to the lack of algorithms to train neural networks efficiently. As a result, they lag behind classical methods in accuracy and efficiency, which prevent most deep neural network-based methods from being not only usable but also efficient.

In recent years, some emerging methods based on shallow neural networks have received increasing attention, such as the extreme learning machine (ELM) \cite{huang2006extreme} and the random feature method (RFM) \cite{2022rfm,2024rfm}. In these methods, all parameters of neural networks except for the output layer are randomly generated and fixed, which simplifies the optimization problem into a linear system. By using the least-squares method to solve the linear system, the expensive training process can be carried out at a low cost with even higher accuracy. These methods and their variants \cite{2021eloc,2023lrnnpg,2024elocDG} have been widely used for solving partial differential equations \cite{2022rfm,MR4664651,li2023local,MR4636189}. However, so far, there is little application of these methods for solving eigenvalue problems.

From another perspective, the methods above fundamentally offer a method to construct basis functions using neural networks. However, as the basis functions constructed by neural networks are sensitive to the choice of hyperparameters, improper choices can significantly impact the accuracy of the approximate solution. Randomly generating these hyperparameters may be inefficient for some problems, which means that the search for more reasonable ways to select and update hyperparameters is a worthwhile pursuit.

To address the impact of hyperparameter selection on solutions, a modified batch intrinsic plasticity method for pre-training the randomly generated hyperparameters in ELM is proposed in \cite{2021eloctrain}. A differential evolution algorithm is employed to compute the optimal or near-optimal parameter $R_m$, where the parameters of the neural networks are initialized with random values generated in the interval $[-R_m,R_m]$ in \cite{2022elocR}. Additionally, a better way to randomly select hyperparameters based on some geometric explanations is provided in \cite{zhang2024transferable}.

There are also some deterministic algorithms to calculate hyperparameters. An alternating iteration process is introduced in \cite{ainsworth2021galerkin,cyr2020robust}, which combines gradient descent for updating the hyperparameters and least squares for solving the linear system. A subspace method based on neural networks (SNN) for partial differential equations  has been proposed in \cite{xu2024subspace}. The main idea of this method is to train neural network-based functions by minimizing a loss function, and then use least squares to find the approximation of the original problem in the subspace spanned by those trained functions, with complete decoupling between training the subspace and solving the linear system. Unlike ELM-like methods and random feature methods, the parameters of the neural network-based functions in SNN are not randomly generated; instead, they are obtained by minimizing some loss function. In comparison to classical deep neural network-based discretization methods, the training process in the SNN focuses on finding better basis functions, which requires relatively lower accuracy in the training process and often needs significantly fewer training epochs. Compared with methods updating hyperparameters by alternating iteration, the updating of hyperparameters and solution of the linear system in SNN is totally decoupled, which greatly reduces the computational cost. Consequently, SNN can achieve high accuracy with low cost in training for solving partial differential equations. Furthermore, a weak form version of SNN has also been proposed in \cite{liu2024subspace}.

In another view, the training process of SNN essentially offers a deterministic approach to update the parameters of the neural network-based functions. The solving process of SNN, which finds an approximate solution in the subspace spanned by the trained functions, can be viewed as a kind of finite dimensional approximation commonly employed in classical methods. The difference lies in the selection of basis functions: classical methods like the finite element method employ piecewise polynomials as basis functions, spectral method employs orthogonal polynomials as basis functions, whereas SNN employs trained neural network-based functions as basis functions.

In this paper, motivated by the idea of using deep learning method to train the basis functions introduced in \cite{xu2024subspace}, we  propose a subspace method based on neural networks for eigenvalue problems. Besides using the deep  learning method to train some neural network-based functions, we further  introduce some dimensionality reduction technique to the trained functions and obtain some neural network-based orthogonal basis. Our method mainly consists of two steps. First, we construct a neural network-based orthogonal basis by some deep learning method and dimensionality reduction technique. Then, we use some classical method, e.g., the Galerkin method, to discretize the eigenvalue problem in the subspace  spanned by the orthogonal basis, and then  obtain the approximate solution.

We apply our method to the solution of some typical eigenvalue problems. The numerical results show that our algorithm can achieve very high accuracy with few training epochs, which demonstrates the advantages of our new algorithm in both computational cost and accuracy.

The subsequent parts of our paper are organized as follows. In Section 2, we provide some preliminary knowledge, including some notation and classical results for a linear elliptic eigenvalue problem, and a brief introduction of the subspace method based on neural networks for some source problems. In Sections 3, we introduce the general framework of our subspace method based on the neural networks for eigenvalue problems and the details for each step,especially on how to construct neural network-based orthogonal basis. In Section 4, we provide some numerical results to show the efficiency of our new method. Finally, we conclude with a brief conclusion.

\section{Preliminaries}
Before introducing our method, we first give some preliminary knowledge, including some basic results of a elliptic eigenvalue problem and a brief introduction to the subspace method based on neural networks for partial differential equations proposed in \cite{xu2024subspace}. 
\subsection{A linear elliptic eigenvalue problem and its finite dimensional discretization}
Let ${\rm \Omega}\subset\mathbb{R}^d(d\geq1)$ be a  polygonal domain. We shall use the standard notation for Sobolev spaces $W^{s,p}(\rm \Omega)$ and their associated norms (see, e.g. \cite{adams2003sobolev}). We denote $H_0^1({\rm\Omega})=\{v\in H^1({\rm\Omega}):v|_{\partial {\rm\Omega}}=0\}$, $(\cdot,\cdot)$ be the standard $L^2$ inner product and $\langle\cdot,\cdot\rangle$ be the standard Euclidean inner product. Additionally, in this paper, we will use the letter $C$ to denote a generic positive constant which may stand for different values at its different occurrences. For convenience, following \cite{xu1992iterative}, the symbols $\lesssim$ are used in this paper, where  $x\lesssim y$ means that $x\leq Cy$ for some constant $C$ that is independent of the spatial discretization.

Consider the following elliptic eigenvalue problem: Find $(\lambda, u) \in \mathbb{R} \times V$ such that $b(u, u)= 1$ and
\begin{equation}
	a(u, v) = \lambda b(u, v), \quad \forall v \in V.\label{AEP}
\end{equation}
Here, $(V,\Vert\cdot\Vert)$ is a Hilbert space with inner product ($\cdot$,$\cdot$) and norm $\Vert\cdot\Vert$. The bilinear forms $a(\cdot,\cdot)$ and $b(\cdot,\cdot)$ on $V\times V$ are symmetric and satisfy $$a(u,v)\lesssim\Vert u\Vert\ \Vert v\Vert,\forall u,v\in V,$$
$$\Vert v\Vert^2\lesssim a(v,v),\forall  v\in V\ \text{and}\ 0<b(v,v), \forall v\in V,v\neq 0,$$

In further, we define the norms $\Vert \cdot\Vert_a$ and $\Vert \cdot\Vert_b$ based on the bilinear forms $a(\cdot,\cdot)$ and $b(\cdot,\cdot)$ as follows:
$$\Vert \cdot\Vert_a=\sqrt{a(\cdot,\cdot)},\Vert \cdot\Vert_b=\sqrt{b(\cdot,\cdot)}.$$
Suppose the norm $\Vert \cdot\Vert_b$  is compact with respect to $\Vert \cdot\Vert$ (\cite{babuvska1989finite}), i.e., from any sequence which is bounded in $\Vert \cdot\Vert$, {one can extract a subsequence which is} Cauchy in $\Vert \cdot\Vert_b$. 

It is well known that the eigenvalue problem (\ref{AEP}) has an eigenvalue sequence $\{\lambda_k\}$ (\cite{BABUSKA1991641})
$$0 < \lambda_1 \leq \lambda_2 \leq \cdots \leq \lambda_k \leq \cdots,\lim_{k \to \infty} \lambda_k = \infty,$$ and corresponding eigenfunctions $$u_1, u_2, \cdots, u_k, \cdots,$$ which satisfy $$b(u_i, u_j) = \delta_{ij},i,j=1,2,\cdots,k,\cdots.$$  In the sequence $\{\lambda_k\}$, $\lambda_k$ are repeated according to their geometric multiplicity.

{
	Let  $V_h$ be a finite dimensional subspace of $V$. The finite dimensional discretization of (\ref{AEP}) in $V_h$ is as follows: Find $(\lambda, u_h) \in \mathbb{R} \times V_h$ such that $b(u_h,u_h)= 1$ and
	\begin{equation}
		a(u_h, v_h) = \lambda b(u_h, v_h), \quad \forall v_h \in V_h.\label{DEP}
	\end{equation}
}
Suppose $V_h=\text{span}\{\varphi_1,\cdots,\varphi_M\}$, where $\{\varphi_i\}_{i=1}^{M}$ is a basis of the finite dimensional subspace $V_h$. Substituting $u_h=\sum_{j=1}^{M}w^{(j)}\varphi_j\ \text{and choosing}\ v_h\ \text{to be}\ \varphi_i$ into (\ref{DEP}), we then have:
$$
\sum\limits_{j=1}^Mw^{(j)}a(\varphi_j,\varphi_i)= \lambda\sum\limits_{j=1}^Mw^{(j)}b(\varphi_j,\varphi_i),i=1,\cdots,M.
$$

Denote $A$ the stiffness matrix and $B$ the mass matrix, that is, $A=(a(\varphi_i,\varphi_j))_{ij}\in\mathbb{R}^{M\times M},B=(b(\varphi_i,\varphi_j))_{ij}\in\mathbb{R}^{M\times M}$. The eigenvalue problem (\ref{AEP}) can be formulated as the following algebraic eigenvalue problem:
\begin{equation}
	A\bm w=\lambda B\bm w,\label{MEP}
\end{equation}
where $\bm w=(w^{(1)},\cdots,w^{(M)})^T\in\mathbb{R}^M.$

By solving the algebraic eigenvalue problem (\ref{MEP}), an approximate solution to the eigenvalue problem (\ref{AEP}) is then obtained.
\subsection{Subspace method based on neural networks}
In this subsection, we will give a brief review for the subspace method based on neural networks proposed in \cite{xu2024subspace} for the partial differential equations, which has been {proven to be able to achieve high accuracy with low cost}.

In \cite{xu2024subspace},  the following general source problem is considered: 
$$
\begin{cases}
	\mathscr{A} u = f & \text{in } \rm\Omega, \\
	\mathscr{B} u = g & \text{on } \partial \rm\Omega,
\end{cases}
$$
where $\mathscr{A}$ and $\mathscr{B}$ are two  differential operators, $f$ and $g$ are given functions.

The neural networks architecture used in \cite{xu2024subspace} is comprised of four parts: an input layer, several hidden layers, a subspace layer and an output layer. For simplicity, we take the case of single output as an example to show the mathematical expression of the neural networks. The expression of multiple outputs can be obtained similarly. 

Denote $L$ the number of hidden layers, $n_0=d$ the number of neurons in the input layer, $n_1,\cdots,n_L$ the number of neurons in each hidden layer, $n_{L+1}=M$ the number of neurons in the subspace layer. The neural networks used in \cite{xu2024subspace} can be expressed as follows:
\begin{equation}\label{NN}
	\left\{
	\begin{aligned}
		& \bm y_0 =   \bm x, \\
		&  \bm y_l =   \sigma \circ \mathcal{T}_l(\bm y_{l-1}),l=1,2,\cdots,L+1, \\
		&	\bm\varphi  =                \bm y_{L+1}, \\
		&u =                    \bm \varphi\cdot \bm w. \\
	\end{aligned}
	\right.
\end{equation}
Here, $\mathcal{T}_l(\bm y) = W_l\cdot \bm y+\bm b_l$, $W_l\in\mathbb{R}^{n_l\times n_{l-1}}$ and $\bm b_l\in\mathbb{R}^{n_l}$ are the weights and biases of the corresponding layer, respectively, $\bm x\in \mathbb{R}^d$ is the input, $\bm\varphi =(\varphi_1,\cdots,\varphi_M)^T$ and $\bm w=(w^{(1)},\cdots,w^{(M)})^T$, $\sigma:\mathbb{R}\rightarrow\mathbb{R}$ is an activation function which can be extended to $\mathbb{R}^n\rightarrow\mathbb{R}^n $ as follows:
$\sigma( \bm z) = [\sigma(z_1),\cdots,\sigma(z_n)]^T \text{for any}\ \bm{z}=(z_1,\cdots,z_n)^T\ \in\mathbb{R}^n.$ 
Denote $\theta=\{W_1,\cdots,W_{L+1},\bm b_1,\cdots,\bm b_{L+1}\}\in {\rm\Theta}:=\mathbb{R}^{n_1\times n_0}\times\cdots\times\mathbb{R}^{n_{L+1}\times n_L}\times\mathbb{R}^{n_1}\times\cdots\times\mathbb{R}^{n_{L+1}}$  the parameters of the neural network-based functions $\{\varphi_i\}_{i=1}^{M}$, and $u(\bm x;\theta,\bm{w})=\sum\limits_{i=1}^M w^{(i)}\varphi_i(\bm x;\theta)$  are the output of neural networks with respect to the input $\bm x$ and parameters $\theta$ and  $\bm w$.

Next, we will review the general framework of the subspace method based on neural networks for the partial differential equations proposed in \cite{xu2024subspace}. Firstly, the neural networks architecture introduced above is adopted, and parameters $\theta$ and $\bm w$ are initialized. Next, train the neural network-based functions, i.e., update parameters $\theta$ by minimizing some loss function, aiming to obtain functions that are better than the initial ones through training. After training, find an approximate solution in the subspace spanned by the trained neural network-based functions $\{\varphi_i\}_{i=1}^{M}$, i.e., update parameters $\bm w$ by solving some linear system, aiming to obtain an approximate solution of the partial differential equation. The above process is summarized as a general framework for solving the partial differential equations in \cite{xu2024subspace}. We also list it here for the  completeness, see Algorithm \ref{alg:A1} for the details.

\begin{algorithm}
	\caption{General framework of the  subspace method based on neural networks for partial differential equations}
	\label{alg:A1}
	\begin{algorithmic} 
		
		\STATE{1. Initialize the neural network architecture, give and then fix  $\bm w_0$.}
		\STATE{2. Define some loss function $\mathcal{L}(\theta,\bm w_0)$ using information contained in the PDEs, then obtain parameters $\theta$ by minimizing the loss function using some optimization method 
			$$\theta^{\ast}=\mathop{\arg \min}\limits_{\theta\in\rm\Theta}\mathcal{L}(\theta,\bm w_0),$$ 
			and obtain the updated neural network-based functions $\{\varphi_i(\bm{x};\theta^{\ast})\}_{i=1}^M$}
		\STATE{3.  Using some classical method to find the approximation for the solution of the original problems in the trained subspace spanned by $\{\varphi_i(\bm{x};\theta^{\ast})\}_{i=1}^M$, i.e., obtain optimal $\bm{w}^{\ast}$,  and then  obtain the approximate solution $u(\bm{x};\theta^{\ast},\bm{w}^{\ast})$.}
		
	\end{algorithmic}
\end{algorithm}

\section{Subspace method based on neural networks for eigenvalue problems}
We consider the following typical eigenvalue problem: 
$$
\begin{cases}
	L u = \lambda u & \text{in } \rm\Omega, \\
	u = 0 & \text{on } \partial \rm\Omega,
\end{cases}
$$
where the differential operator $L:=-\nabla\cdot(\alpha\nabla)+\beta$ and $\alpha=(\alpha_{ij})_{ij}$ is symmetric positive definite with $\alpha_{ij}\in W^{1,\infty}({\rm\Omega}) (i,j=1,\cdots,d)$, and $0\leq\beta\in L^{\infty}(\rm\Omega).$ 

The associated weak form of the eigenvalue problem above is as follows: Find $(\lambda, u) \in \mathbb{R} \times H_0^1(\rm\Omega)$ such that $b(u, u)=1$ and
\begin{equation}\label{EP}
	a(u, v) = \lambda b(u, v), \quad \forall v \in H_0^1(\rm\Omega),
\end{equation}
which is nothing but (\ref{AEP}) with  $a(u,v)=(\alpha \nabla u,\nabla v)+(\beta u,v)$ and $b(u,v)=(u,v)$. Here, $a(\cdot,\cdot)$ and $b(\cdot,\cdot)$ are two positive definite symmetric bilinear forms on $H_0^1({\rm\Omega})\times H_0^1(\rm\Omega)$.

\begin{remark}
	It is worth mentioning that our method is also valid for a more general bilinear form $a(\cdot,\cdot)$ that
	$$
	\Vert v\Vert_1^2-C\Vert v\Vert_0^2\lesssim a(v,v),\quad \forall v\in H_0^1(\rm\Omega)
	$$
	holds for some constant $C$. (see, e.g., Remark 2.9 in \cite{dai2008convergence})
\end{remark}

Let  $V_h$ be a finite dimensional subspace of $V$. The finite dimensional discretization of (\ref{EP}) is as follows: Find $(\lambda, u_h) \in \mathbb{R} \times V_h$ such that $b(u_h,u_h)= 1$ and
\begin{equation}
	a(u_h, v_h) = \lambda b(u_h, v_h), \quad \forall v_h \in V_h.\label{DEP1}
\end{equation}

In this paper, we consider to obtain the smallest $k$ eigenvalues and their corresponding eigenfunctions of the eigenvalue problem (\ref{EP}) and discretize the eigenvalue problem (\ref{EP}) in a finite dimensional space spanned by a neural network-based orthogonal basis. 
Our basic idea is as follows: we 
first use deep learning method  to train some neural network-based functions $\{\varphi_i\}_{i=1}^M$ which is similar to \cite{xu2024subspace}; then, we introduce some dimensionality reduction technique to obtain some orthogonal basis $\{\psi_i\}_{i=1}^K$ from these trained neural network-based functions $\{\varphi_i\}_{i=1}^M$; after that, we use some classical method to obtain the approximation solution of the eigenvalue problem (\ref{EP}) in the subspace $V_h=\text{span}\{\psi_1,\cdots,\psi_K\}$.  Usually, there holds $K\ll M$. 

Here, we would like to point out that different from the work in \cite{xu2024subspace}, here, we add a step using some dimensionality reduction technique to the trained neural network-based functions. In fact, we have done some numerical tests, which show the  necessity of using dimensionality reduction technique. Some of the results are provided in Appendix B.

We first give a general framework of our subspace method based on neural networks for solving eigenvalue problems, as in Algorithm \ref{alg:A2},   and then provide more details for each of the steps in the following parts of this section. Here, to make the framework more readable, we combine the training of the neural network-based functions using some deep learning method and obtaining the orthogonal basis from those trained function by some dimensionality reduction technique together and view them as a separate step for obtaining an orthogonal basis. 
\begin{algorithm}
	\caption{General framework of subspace method based on neural networks for eigenvalue problems}
	\label{alg:A2}
	\begin{algorithmic}
		
		\STATE{1. Construct a neural network-based orthogonal basis $\{\psi_i\}_{i=1}^{K}$ by some deep learning method and dimensionality reduction technique.}
		\STATE{2. Use some classical method, e.g., the Galerkin method, to discretize the  eigenvalue problem (\ref{EP}) in the subspace spanned by the basis $\{\psi_i\}_{i=1}^K$ and  obtain an approximate solution.}
	\end{algorithmic}
\end{algorithm}
\subsection{Constructing a neural network-based orthogonal basis}
In this subsection, we  introduce how we obtain the neural network-based orthogonal basis $\{\psi_i\}_{i=1}^{K}$ by using some deep learning method and dimensionality reduction technique. 

Inspired by the subspace method based on neural networks proposed in \cite{xu2024subspace}, we consider to use neural network-based functions to construct basis and train the parameters of those neural network-based functions by minimizing some loss function. Furthermore, we observe the difference between solving eigenvalue problems and source problems. To solve the algebraic eigenvalue problem resulted from the operator eigenvalue problems with high accuracy, the stiffness matrix and the mass matrix generated by neural network-based basis functions are required to be symmetric and have good condition numbers. Usually, the trained neural network-based functions lack good linear independence. Consequently, some strategy must be taken to remove the redundant information contained in the trained neural network-based functions. In fact, as referred to in the beginning of this section, we have done some tests to show the necessity of using dimensionality reduction technique in Appendix B.  In addition, to deal with boundary conditions, we also hope the basis satisfies boundary conditions of eigenvalue problem (\ref{EP}) exactly. 

Based on the consideration above, we get the main steps for obtaining a neural network-based orthogonal basis $\{\psi_i\}_{i=1}^{K}$ by using some deep learning method and dimensionality reduction technique as follows: First, choose a neural networks architecture with its parameters to be determined; then, modify neural network-based functions to make them satisfy boundary conditions exactly; after that, construct a proper loss function by the neural network-based functions and the information provided by the problem to be solved, and obtain the undetermined parameters in the neural network-based functions by minimizing the loss function; finally, construct a neural network-based orthogonal basis from the trained neural network-based functions by using some dimensionality reduction technique. In the following of this subsection, we will introduce the details for each step.

\subsubsection{\bf Neural networks architecture} 

The neural networks architecture used here is similar to the one introduced in the previous section. The only difference is that the neural networks in Section 2 is used to approximate the solution to the  source problem, here the neural networks are used to approximate the eigenfunctions corresponding to the first $k$ eigenvalues of some partial differential operator. Therefore, the dimension of the output layer is $k$. Thus, the output of {the} neural networks can be presented as follows,
$$\bm u(\bm x;\theta, w) = (\bm\varphi\cdot \bm w_1,\cdots,\bm\varphi\cdot \bm w_k)^T = \left(\sum\limits_{j=1}^Mw_1^{(j)}\varphi_j(\bm x,\theta),\cdots,\sum\limits_{j=1}^Mw_k^{(j)}\varphi_j(\bm x,\theta)\right)^T,
$$
where $ w=(\bm w_1,\cdots,\bm w_k)$, $\bm w_s=(w_s^{(1)},\cdots,w_s^{(M)})^T $ and $u_s(\bm{x}; \theta, \bm w_s)$ denotes the $s$-th output of the neural networks for $s=1,\cdots,k$. Similar to those introduced in Section 2, here,  $\bm \varphi(\bm x,\theta)=\sigma\circ\mathcal{T}_{L+1}\circ\cdots\circ\sigma\circ\mathcal{T}_1(\bm x)$, {$\sigma$ is an activation function, $\mathcal{T}_l(\bm y) = W_l\cdot \bm y+\bm b_l$} and $\theta=\{W_1,\cdots,W_{L+1},\bm b_1,\cdots,\bm b_{L+1}\}$ are the parameters of the neural network-based functions $\{\varphi_i\}_{i=1}^{M}$.

The initialization of the parameters $\theta$ and  $\bm w_1,\cdots,\bm w_k$ are as follows. For the parameters $\theta$, we adopt the default initialization strategy provided by PyTorch. For the parameters $\bm w_1,\cdots,\bm w_k$, there are various methods for assigning initial values. Numerical experiments provided in \cite{xu2024subspace} indicate that different choices of $\bm w_1,\cdots,\bm w_k$ have no obvious impact on the results. Therefore, here, we do not pay much attention on the choice of $\bm w_1,\cdots,\bm w_k$. We just give some simple choices. That is, for a single output, we choose $\bm w=(1,\cdots,1)^T\in\mathbb{R}^M$. For multiple outputs, a random approach is chosen, that is, $\bm w_1,\cdots,\bm w_k\in\mathbb{R}^M$ are assigned random values generated from a uniform distribution defined over the interval $[-1,1]$.

\subsubsection{Handling the boundary condition}
For eigenvalue problems, some boundary conditions are usually required. Since the neural network-based functions are usually nonlocal, it is not so easy to deal with the boundary condition by using the same way as classical methods. Here, we choose to modify those functions to make them satisfy the boundary condition. For example,  for the case of zero boundary condition with  $\rm\Omega$ being bounded, we  can modify the basis function by the following way. We multiply $\varphi_i$ by a function $f$ such that 
\begin{equation}
	\bar{\varphi}_i = \varphi_i f,\ \bar{\varphi}_i(\bm x)=0,\  \forall \bm x\in \partial{\rm\Omega},\  \forall i=1,\cdots,M.\label{M1}
\end{equation}
For example, for the case of the 2D unit square domain, we can choose $f(\bm x) = x_1x_2(1-x_1)(1-x_2)$. 

If $\rm\Omega$ is unbounded, it is necessary to modify the neural network-based functions to ensure that the modified functions are $L^2$ integrable. We can achieve this by multiplying $\varphi_i$ by a function $f$ such that 
\begin{equation}
	\bar{\varphi}_i = \varphi_i f,\ \int_{\rm\Omega}\bar{\varphi}_i^2dx<\infty,\ \forall i=1,\cdots,M.\label{M2}
\end{equation}
Since the Hermite-Gauss quadrature scheme is used to numerically approximate the {integrals on an unbounded domain} in this paper, we can choose $f(\bm x) = \exp(-\frac{1}{2} \bm x^T\bm x)$.

For simplicity, in the remaining parts of this paper, $\varphi_i$ will refer to the modified neural network-based function $\bar{\varphi}_i$ as described above.

\subsubsection{The loss function and the training}
In this subsection, we will introduce our choice for the loss function $\mathcal{L}$ and a specific approach to train the  neural network-based functions.

We first introduce our choices for the loss function. Before that, we first introduce some properties of the eigenvalue problem (\ref{EP}) to give us an idea on how to choose the loss function. 

Denote $V_k = \text{span}\{v_1, \cdots, v_k\}$ an arbitrary $k$-dimensional subspace of $V$ and $U_k = \text{span}\{u_1, \cdots, u_k\}$ the eigensubspace corresponding to the smallest $k$ eigenvalues.

According to the minimum-maximum principle \cite{BABUSKA1991641}, the sum of the smallest $k$ eigenvalues and eigenfunctions corresponding to the smallest $k$ eigenvalues satisfy the following equations,
\begin{equation}\label{EV}
	\sum_{i=1}^k \lambda_i = \min_{V_k=\text{span}\{v_1, \cdots, v_k\} \atop v_k \in V, dim(V_k)=k}  \text{Trace} \left( \mathcal{B}^{-1}(v_1, \cdots, v_k)\mathcal{A}(v_1, \dots, v_k) \right),
\end{equation}
\begin{equation}\label{EF}
	\left( u_1,\cdots,u_k\right)=\mathop{\arg \min}_{V_k=\text{span}\{v_1, \cdots, v_k\} \atop v_k \in V, dim(V_k)=k}  \text{Trace} \left( \mathcal{B}^{-1}(v_1, \cdots, v_k)\mathcal{A}(v_1, \dots, v_k) \right). 
\end{equation}
Here, $\mathcal{A}$ is the stiffness matrix and $\mathcal{B}$ is the mass matrix, that is, $\mathcal{A}(v_1, \dots, v_k)=(a(v_i,v_j))_{ij}\in\mathbb{R}^{k\times k}$, $\mathcal{B}(v_1, \dots, v_k)=(b(v_i,v_j))_{ij}\in\mathbb{R}^{k\times k}$.

We now consider the approximation of the eigenvalue problem (\ref{EP}) on the finite dimensional subspace $V_h$ spanned by neural network-based functions. Set $V_h = \text{span}\{\varphi_1(\bm x;\theta),\cdots,\varphi_M(\bm x;\theta)\}\subset H^1_0(\rm\Omega)$.	Set $v_i(\bm x;\theta,\bm w_i)=\sum\limits_{j=1}^M w_i^{(j)}\varphi_j(\bm x;\theta)$, $i=1, \cdots, k$. 	Motivated by (\ref{EV}) and (\ref{EF}), we then construct the following loss function for finding first $k$ eigenpairs (see also in \cite{WANG2024112928}):
\begin{equation}
	\mathcal{L}(\theta,w)=\text{Trace}(\mathcal{B}^{-1}(v_{1}(x;\theta,\bm w_1),\cdots,v_{k}(x;\theta,\bm w_k))\mathcal{A}(v_{1}(x;\theta,w_1),\cdots,v_{k}(x;\theta,\bm w_k))).\label{loss}
\end{equation} 

Now, we turn to introduce how we train the neural network-based functions. 
Similar to Algorithm \ref{alg:A1} proposed in  \cite{xu2024subspace}, for training the neural network-based functions, we first fix the parameter  $w$, that is, set $w =  w^{(0)}$ with $w^{(0)}$ being given at initialization and then obtain  parameters $\theta$  by solving the following minimization problem:
\begin{equation}
	\min_{\theta\in\rm\Theta}\mathcal{L}(\theta, w^{(0)}).\label{O}
\end{equation}

The gradient descent method is adopted to solve the minimization problem (\ref{O}). That is, we update the parameters $\theta$ by 
$$ \theta^{(l+1)}=\theta^{(l)} - \eta \nabla_{\theta} \mathcal{L}(\theta^{(l)},  w^{(0)}), $$ 
where $\theta^{(l)}$ is the value of $\theta$ in $l$-th epoch and $\eta$ is the learning rate adjusted by the Adam optimizer \cite{kingma2014adam}. Here, when minimizing the loss function, a fixed $w^{(0)}$ is adopted

We have also briefly explained how the training process affects the accuracy of final solution through some numerical experiments in Appendix A.

For the choice of termination criteria of the gradient descent method, we notice that the goal of training is to obtain the basis other than the solution, there is no need to train until the loss is small enough, e.g. find the $\theta^{\ast}$ s.t. $\mathcal{L}(\theta^{\ast}, w^{(0)})-\min\limits_{\theta\in\rm \Theta}\mathcal{L}(\theta,w^{(0)})\leq\epsilon$. Instead, in \cite{xu2024subspace}, the following termination criterion is adopted:
$$\frac{\mathcal{L}(\theta, w^{(0)})}{\mathcal{L}_0(\theta, w^{(0)})}\leq \epsilon\ \text{or}\ N_{epoch}\geq N_{max},$$
where $\mathcal{L}_0(\theta,w^{(0)})$ is the initial loss, $\epsilon$ is the tolerance, $N_{epoch}$ and $N_{max}$ are the current and maximum number of training epochs, respectively. However, since the minimum value of the loss function we used in this paper is not 0, this criterion may not be suitable for determining when to stop training. As an alternative, we propose a new termination criterion.

We note that for an iterative method, the error between two adjacent iterations is usually used to determine when to terminate the iteration. That is, the iteration is terminated when the error between two adjacent iterations is less than some given tolerance $\epsilon$. Note that the loss usually decreases in an oscillatory manner during the training, hence using the change of the loss function during multiple adjacent epochs as an indicator may be better than using the change during two adjacent epochs. Therefore, we choose to use the change of the loss during multiple adjacent epochs as termination criterion. In this paper, we adopt the change of the loss during $m$ adjacent epochs as an indicator. Any other reasonable choices are also acceptable.

Let  $l_s =\mathcal{L}_s(\theta, w)$ be the loss in the $s$-th epoch, $\bar l_s =\sum\limits_{i=1}^{m} \vert l_{s-i}$ $-l_{s+1-i}\vert$ be the change of the loss function during $m$ adjacent epochs. The following termination criterion is adopted:
\begin{equation}\label{TC}
	\vert\frac{\bar l_s}{l_s}\vert\leq  \epsilon, \ \text{or}\ N_{epoch}\geq N_{max},
\end{equation}
where $\epsilon$ is a given tolerance, $N_{epoch}$ and $N_{max}$ are  the current and maximum training epochs, respectively. Under this criterion, the termination of training means that the relative change of the loss function during $m$ adjacent epochs has been relatively small, which indicates that the benefits of continued training may be small. Since there is no need to train until the loss is small enough, a relatively big $\epsilon$ and a relatively small $N_{max}$ can be chosen.

\subsubsection{Dimensionality reduction technique} 
In this subsection, we will introduce how we obtain the neural network-based orthogonal basis $\{\psi_i\}_{i=1}^{K}$ from the trained neural networks-based functions $\{\varphi_i\}_{i=1}^M$ which may numerically lack good linear independence while keeping as much necessary information as possible by some dimensionality reduction technique. Here, we choose the proper orthogonality decomposition (POD) to achieve this goal. 


We can describe the problem as follows. Denote $V=\text{span}\{\varphi_1,\cdots,\varphi_M\}$, $V_k \subset V$ any $K$-dimensional subspace of $V$, and ${\rm\Pi}_K$ the orthogonal projection operator from $V$ to $V_K$. The problem is to find an optimal $K$-dimensional subspace $V_K^{\ast}$ to approximate the original space $V$ good enough. There will be two subproblems. One is that, for a given dimension of the subspace $K\leq M$, how to find an optimal $K$-dimensional subspace $V_K^{\ast}$ to approximate the original space $V$ in a least–squares sense, the other is how to choose a suitable dimension of the subspace $K$. For the first problem, i.e. we seek an optimal orthogonal projection ${\rm\Pi}_K^{\ast}:V\rightarrow V_K^{\ast}$ minimizing the projection error of ${\rm \Phi}=(\varphi_1,\cdots,\varphi_M)^T\in[L^2({\rm\Omega})]^M$ under the following least–squares distance,
$$
\Vert{\rm\Phi}-{\rm\Pi}_K{\rm\Phi}\Vert^2:=\sum\limits_{i=1}^M \big|\big|\varphi_i-\sum\limits_{j=1}^K(\varphi_i,\psi_j)\psi_j\big|\big|_{L^2({\rm\Omega})}^2.
$$

Hence, we need to solve the following minimization problem,
\begin{equation}\label{OPL2}
	\min_{\psi_1,\cdots,\psi_K\in {L^2}({\rm\Omega})}\sum\limits_{i=1}^M \big|\big|\varphi_i-\sum\limits_{j=1}^K(\varphi_i,\psi_j)\psi_j\big|\big|_{L^2({\rm\Omega})}^2\ \text{s.t.}\ (\psi_i,\psi_j)=\delta_{ij},
\end{equation}

Denote $U=[u_1,\cdots,u_K]\in\mathbb{R}^{M\times K}$ and 
\begin{equation}
	{\rm\Psi}=U^T{\rm\Phi}.\label{Psi}
\end{equation}
We can rewrite the minimization problem (\ref{OPL2}) to the minimization problem (\ref{OPR}) through the following simple derivation,
$$
\sum\limits_{i=1}^M \big|\big|\varphi_i-\sum\limits_{j=1}^K(\varphi_i,\psi_j)\psi_j\big|\big|_{L^2({\rm\Omega})}^2=\sum\limits_{i=1}^M \big|\big|e_i-\sum\limits_{j=1}^K\langle e_i,u_j\rangle _{\mathcal{M}}u_j\big|\big|_{\mathcal{M}}^2,
$$
where $\mathcal{M}=((\varphi_i,\varphi_j))_{i,j}\in\mathbb{R}^{M\times M}$, $\langle x,y\rangle_{\mathcal{M}}=\langle x,\mathcal{M}y\rangle,\Vert x\Vert_{\mathcal{M}}=\sqrt{\langle x,x\rangle_{\mathcal{M}}}, \forall x,y\in\mathbb{R}^M$ and $e_1,\cdots,e_M$ are the standard basis in $\mathbb{R}^M$.

\begin{equation}\label{OPR}
	\min_{u_1,\cdots,u_K\in \mathbb{R}^M}\sum\limits_{i=1}^M \big|\big|e_i-\sum\limits_{j=1}^K\langle e_i,u_j\rangle_{\mathcal{M}} u_j\big|\big|_{\mathcal{M}}^2\ \text{s.t.}\ \langle u_i,u_j\rangle_{\mathcal{M}}=\delta_{ij},
\end{equation}

Notice that the problem (\ref{OPR}) is nothing but the minimization problem considered in the POD method, which has been systematically studied (see \cite{pinnau2008model,volkwein2011model}), hence we give the following optimal orthogonal projection ${\rm\Pi}_K^{\ast}$ of the problem (\ref{OPR}) directly,
\begin{equation}\label{POD1}
	\mathcal{M}^2u_j=\mu_j\mathcal{M}u_j,j=1,\cdots,K,
\end{equation}
\begin{equation}\label{POD2}
	\min_{u_1,\cdots,u_K\in \mathbb{R}^M}\sum\limits_{i=1}^M \big|\big|e_i-\sum\limits_{j=1}^K\langle e_i,u_j\rangle_{\mathcal{M}} u_j\big|\big|_{\mathcal{M}}^2=\sum\limits_{j=K+1}^M\mu_j,
\end{equation}
where $\mu_1\geq\cdots\mu_M\geq 0$ are the first $M$ eigenvalues of $\mathcal{M}^2$ relative to $\mathcal{M}$, $u_1,\cdots,u_M$ are the corresponding eigenvectors.

The vectors {$u_1,\cdots,u_K$} are called POD modes. In practical implementation, we can obtain them by solving the following symmetric $M\times M$ eigenvalue problem$$
\mathcal{M} \bar u_j=\mu_j\bar u_j,j=1,\cdots,K,$$
$$u_j=\frac{1}{\sqrt{\mu_j}}\bar u_j,j=1,\cdots,K.$$

\begin{remark}
	In actual implementation, the condition number of $\mathcal{M}$ is usually very large, and it is difficult to obtain a high-precision solution by directly solving the eigenvalue problem. Instead, we choose to calculate the singular value of $\mathcal{M}^{\frac{1}{2}}$, from which we can also get $\mu_j, u_j, j= 1,\cdots,K.$ By the way, we only need to deal with a problem whose condition number becomes the square root of the original problem.
	
	In addition, we do not need to pay any extra cost for calculating $\mathcal{M}^{\frac{1}{2}}$, because it has a natural approximation as follows. Denote $\Phi_h = (\varphi_i(\bm x^{(j)}))_{ij}\in \mathbb{R}^{M\times n}, W_h = diag(\rho_1,\cdots,\rho_n)\in \mathbb{R}^{n\times n}$, where $\{\bm x^{(i)},\rho_i\}_{i=1}^n$ represents the numerical integration point and weight, respectively. For $\mathcal{M}$, we have the following approximation, $\mathcal{M}\approx \Phi_h W_h \Phi_h^T$. Hence, we can use $W_h^{\frac{1}{2}}\Phi_h^T$ to approximate $\mathcal{M}^{\frac{1}{2}}$.
\end{remark}

Regarding the choice of the subspace dimension $K$, since (\ref{POD2}) shows the approximate error of the optimal $K$-dimensional subspace is related to the sum of the last $M-K$ eigenvalues of the matrix $\mathcal{M}$ (or the square root of the eigenvalues), it's reasonable to choose a small enough $K$ such that the first $K$ eigenvalues (or the square root of the eigenvalues) account for a high proportion in the overall, i.e., 
$$
K=\mathop{\arg\min}\{I(N):I(N)\geq1-\gamma\},
$$
where $
I(N)=\frac{\sum\limits_{j=1}^{N}\mu_j }{\sum\limits_{j=1}^{M}\mu_j}
$
\text{or}
$
I(N)=\frac{\sum\limits_{j=1}^{N}\sqrt{\mu_j }}{\sum\limits_{j=1}^{M}\sqrt{\mu_j}}
$
and $\gamma$ is a given parameter.

In our experiment, we find that the indicator $I(N)$ defined by the latter works better.

To see the necessity of employing dimensionality reduction technique, we have done some numerical tests, whose results are present in Appendix B.

\subsubsection{Algorithm} 

By summarizing the steps above, we obtain our algorithm for constructing the neural network-based orthogonal basis by deep learning method and dimensionality reduction technique, see Algorithm \ref{alg:A3} for the details. 
\begin{algorithm}
	\caption{Constructing the neural network-based orthogonal basis}
	\label{alg:A3}
	\begin{algorithmic}
		\renewcommand{\algorithmicrequire}{\textbf{Input:}}
		\renewcommand{\algorithmicensure}{\textbf{Output:}}
		\REQUIRE {the neural network structure: $L, n_0, n_1, \cdots, n_L, M$; parameters used for training the basis:  $\epsilon,m,N_{\max},\gamma$; number of eigenpairs needed: $k$}
		\ENSURE { the number of orthogonal basis functions $K$ and  the basis $\{\psi_i(\bm x,\theta^{\ast})\}_{i=1}^K.$}
		\STATE{1. Choose a neural network architecture and modify neural network-based functions to make it satisfy the boundary condition.}
		\STATE{2. Train the  neural network-based functions $\{\varphi_i(\bm x,\theta)\}_{i=1}^M$ by minimizing the loss function (\ref{loss}) until (\ref{TC}) is satisfied.}
		\STATE{3. Obtain the neural network-based orthogonal basis $\{\psi_i(\bm x,\theta^{\ast})\}_{i=1}^K$ from the trained neural network-based functions $\{\varphi_i(\bm x,\theta^{\ast})\}_{i=1}^M$ by POD method.} 
		
	\end{algorithmic}
\end{algorithm}

\subsection{Galerkin discretization in the neural network-based subspace}
In this subsection, we will introduce the details on how we discretize the  eigenvalue problem (\ref{EP}) in the subspace spanned by the neural network-based orthogonal basis  $\{\psi_i(\bm x,\theta^{\ast})\}_{i=1}^K$. 

Similar to the introduction in Section 2, we consider the following finite dimensional discretization of (\ref{EP}): Find $(\lambda, u_h) \in \mathbb{R} \times V_h$ such that $b(u_h,u_h)= 1$ and
\begin{equation}
	a(u_h, v_h) = \lambda b(u_h, v_h), \quad \forall v_h \in V_h.\label{evp-dis}
\end{equation}

Let 
$V_h = \text{span}\{\psi_1(\bm x;\theta^{\ast}),\cdots,\psi_K(\bm x;\theta^{\ast})\}$.  Set $u_{NN}(\bm x;\theta^{\ast},\bm w)=\sum\limits_{j=1}^K w^{(j)}\psi_j(\bm x;\theta^{\ast})$, and substitute  $u_{NN}(\bm x;\theta^{\ast},\bm w)$ into (\ref{evp-dis}) and choose  $ v_h(\bm x)=\psi_i(\bm x;\theta^{\ast})$, we get 
\begin{eqnarray}\label{evp-dis-snn}
	\sum\limits_{j=1}^K w^{(j)}a(\psi_j,\psi_i)=\lambda\sum\limits_{j=1}^K w^{(j)}b(\psi_j,\psi_i),i=1,\cdots,K.
\end{eqnarray} 

If we set $A=(a(\psi_i(\bm x;\theta^{\ast}),\psi_j(\bm x;\theta^{\ast})))_{ij},B=((\psi_i(\bm x;\theta^{\ast}),\psi_j(\bm x;\theta^{\ast})))_{ij}\in\mathbb{R}^{K\times K},$ $\bm w=(w^{(1)},\cdots,w^{(K)})^T\in\mathbb{R}^K$, then (\ref{evp-dis-snn}) can be written as the following algebraic eigenvalue problem
\begin{equation}\label{MEP1}
	A  \bm w=\lambda B \bm w.
\end{equation}

To calculate the stiffness matrix $A$ and the mass matrix $B$, we have to do some numerical integration. For simplicity, we take $a(u,v)=(\nabla u,\nabla v),b(u,v)=(u,v)$ as an example to show how to compute each element of the stiffness matrix and mass matrix by numerical integration. Since $\theta^{\ast}$ is fixed, for simplicity, we will omit $\theta^{\ast}$ and denote $\psi_i(\bm x,\theta^{\ast})$ as $\psi_i(\bm x)$.They are calculated as follows:
$$
(\psi_i,\psi_j)\approx\sum\limits_{l=1}^n \rho_l \psi_i(\bm x^{(l)})\psi_j(\bm x^{(l)}),
(\nabla\psi_i,\nabla\psi_j)\approx\sum\limits_{l=1}^n \rho_l \nabla\psi_i(\bm x^{(l)})\nabla\psi_j(\bm x^{(l)}).
$$
Here, $\{\bm x^{(l)}\}_{l=1}^n$ are the numerical integration points, $\{\rho_l\}_{l=1}^n$ are the weights for each numerical integration point and $n$ is the number of numerical integration points.

The smallest $k$ eigenvalues and their corresponding eigenfunctions $(\lambda_{h,1},\bm w_1),\cdots,$ $(\lambda_{h,k}, \bm w_k)$ can be obtained by solving the algebraic eigenvalue problem (\ref{MEP1}) using some eigensolver \cite{bai2000templates,saad2011numerical}, then the approximation of the $k$ eigenpairs $(\lambda_{h,1},u_{h,1}),\cdots,$ $(\lambda_{h,k},u_{h,k})$ are obtained by setting $u_{h,i}=\sum\limits_{j=1}^M w_i^{(j)}\psi_j(x,\theta^{\ast}), i=1,\cdots,k.$

\section{Numerical experiments}
In this section, we will use four typical eigenvalue problems, the Laplace eigenvalue problem, the decoupled harmonic oscillator, the {coupled} harmonic oscillator and the hydrogen atom to demonstrate the performance of the methods proposed in our paper. 

The algorithms are implemented based on the deep learning framework PyTorch. All the numerical experiments are carried out on the LSSC-IV platform at the Academy of Mathematics and Systems Science, Chinese Academy of Sciences. 

The following parameter settings are used in our numerical experiments: The number of eigenpairs $k$ is set to 15 in all examples. The neural networks architecture is implemented by a fully connected neural network with 3 hidden layers and a subspace layer with each hidden layer containing 100 neurons and the subspace layer containing $M$ neurons. We choose the activation function $\sigma(x)=\sin(x).$ During the initialization part, the default initialization strategy is adopted to initialize parameters $\theta$ and the random values generated according to a uniform distribution in the range [-1, 1] are applied to initialize parameters $\bm w_1,\cdots,\bm w_k$. The Adam optimizer with the default learning rate is used for training. In all experiments, we adopt the termination criteria introduced in (\ref{TC}) with the parameter $N_{max} = 5000$ and $m=10$. Different $\epsilon$ are chosen for different problems to achieve appropriate training (i.e., fewer epochs while maintaining accuracy). The solver $scipy.linalg.eigh$ in Python is used for solving the algebraic eigenvalue problem (\ref{MEP1}). For reproducibility, all the random seeds used in the numerical experiment are set to 1.

For simplicity, we denote our subspace method based on neural networks for eigenvalue problems, Algorithm 2,  as SNN-EP. To show the advantage of our new method to existing neural network-based methods, we compare our method with other neural network-based method, the one proposed in \cite{WANG2024112928}, which can be viewed as a variant of DRM capable of solving multiple eigenpairs simultaneously, and denote it as DRM-M. 
For the fairness of the comparison, we set the same neural networks architecture, initialization method, and the optimizer. The epochs in SNN-EP are determined by the termination criterion (\ref{TC}), while those in DRM-M are set to $25000$. Setting the epochs big enough to do the training  is a common usage when using deep learning method to solve PDEs. 

The relative error of the approximate eigenvalues and eigenfunctions are defined as follows:
$$
err_{\lambda,l}:=\frac{\vert\lambda_l - \lambda_{h,l}\vert}{\vert\lambda_l\vert},err_{L^2,l}:=\frac{\Vert u_l - u_{h,l}\Vert_{L^2({\rm\Omega})}}{\Vert u_l\Vert_{L^2({\rm\Omega})}},err_{H^1,l}:=\frac{\Vert u_l - u_{h,l}\Vert_{H^1({\rm\Omega)}}}{\Vert u_l\Vert_{H^1({\rm\Omega})}},
$$
where $\lambda_l$ and $u_l$ are reference eigenvalues and eigenfunctions, $l=1,\cdots,k$.
\subsection{Laplace eigenvalue problem}
We first consider the following 2D Laplace eigenvalue problem on a square domain: Find $(\lambda,u)\in \mathbb{R} \times H_0^1(\rm\Omega)$ such that $\int_{\rm\Omega}\vert u\vert^2 d\rm\Omega=1$  and
$$
\begin{cases}
	- {\rm\Delta} u = \lambda u & \text{in } \rm\Omega, \\
	u = 0 & \text{on } \partial \rm\Omega,
\end{cases}
$$
where  ${\rm\Omega} = [0,1]^2$. The exact eigenvalues are as follows:
$$
\lambda_{n_1,n_2} = \pi(n_1^2 + n_2^2) ,n_1, n_2 = 1, 2, \cdots,
$$
and the corresponding eigenfunctions  are
$$
u_{n_1,n_2}(x_1, x_2) = \sqrt{2} \sin \left( \pi n_1 x_1 \right) \sin \left( \pi n_2 x_2 \right), n_1, n_2 = 1, 2, \cdots.
$$

For using SNN-EP to solve the above 2D Laplace eigenvalue problem, since  $\rm\Omega$ is bounded, we choose the Legendre-Gauss integration scheme with 32 integration points in each dimension. For training, we set $\epsilon=10^{-3}$. And for the dimensionality reduction part, we set $\gamma=10^{-11}$.

The dimension of the reduced subspace and the relative errors of the smallest 15 eigenvalues obtained by SNN-EP  with different  $M$ are presented in Figure 1. It can be observed that when   $M$ is small, the dimension of the reduced subspace $K$ increases significantly with the increase of $M$. At this stage, the relative errors of the smallest 15 eigenvalues rapidly decay to below $10^{-12}$. When $M$ reaches $300$, the increase of   $K$ slows down significantly and the relative errors of the smallest 15 eigenvalues remain almost unchanged. Therefore, we only show results obtained by setting $M=300$ for both SNN-EP and DRM-M.

\begin{figure}[htbp]
	\label{fig:laplace_fig}
	\includegraphics[height=7.2cm,width=11.6505cm]{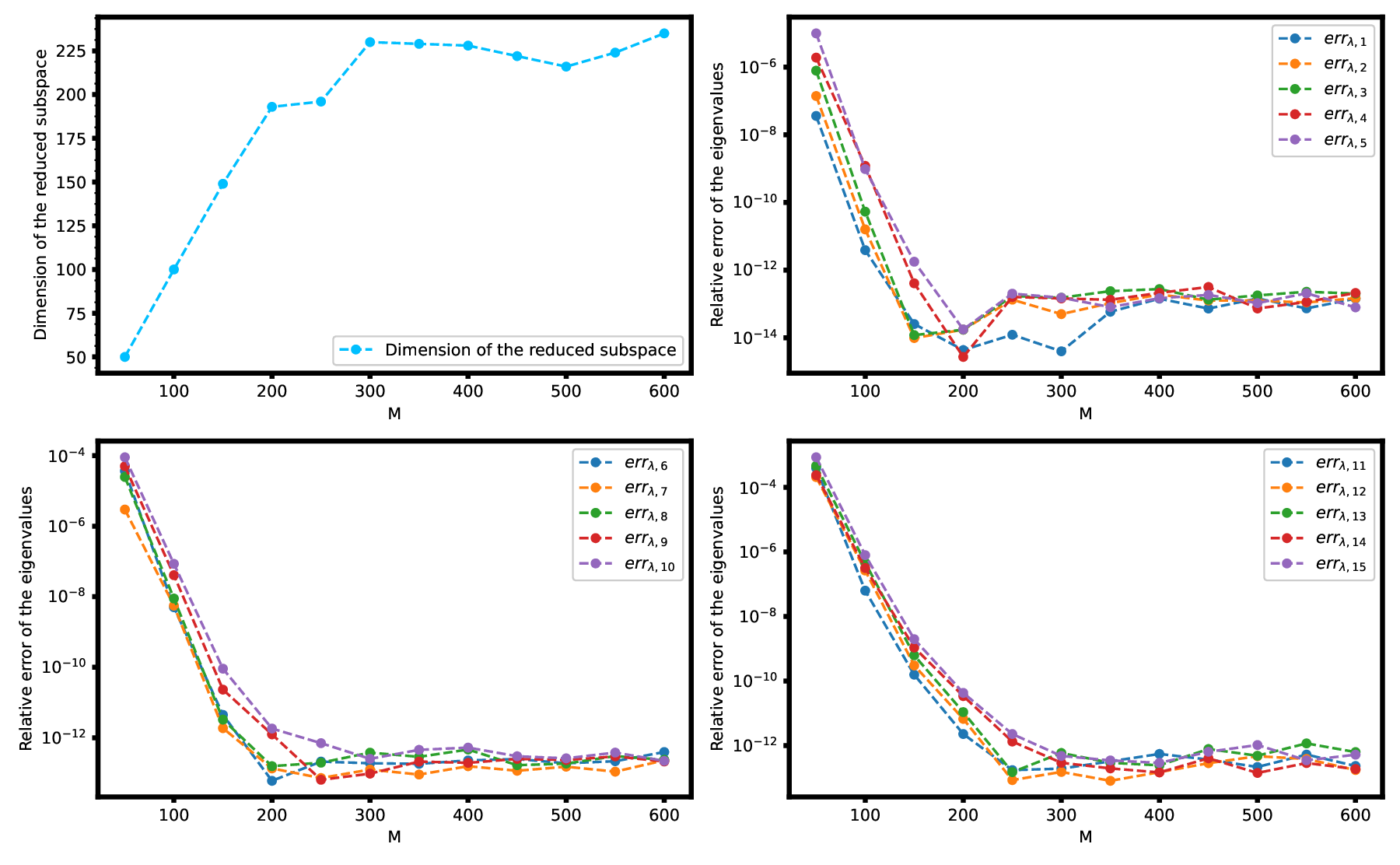}
	\centering
	\caption{The dimension of the reduced subspace and the relative errors of the smallest 15 eigenvalues for 2D Laplace eigenvalue problem obtained by SNN-EP with different $M$.}
\end{figure}

The relative errors of the smallest $15$ eigenpairs obtained using SNN-EP and DRM-M with $M=300$ are presented in Table \ref{tab:laplace_table}, and epochs required for the training process of SNN-EP are shown in Table \ref{tab:laplace_table_epoch}. It can be seen that the accuracy  obtained by DRM-M reaches only $10^{-8}$ for eigenvalues and $10^{-4}$ for eigenfunctions after $25000$ epochs, while our method SNN-EP achieves accuracy below $10^{-12}$ for eigenvalues and nearly $10^{-7}$ for eigenfunctions, only using $100$ to $400$ epochs.
\begin{table}[htbp]
	\centering
	\caption{The relative errors of the smallest 15 eigenvalues  for the 2D Laplace eigenvalue problem obtained by setting $M=300$.}
	\label{tab:laplace_table}
	\scalebox{0.75}{
		\begin{tabular}{ll|lll|lll}
			\hline
			\multicolumn{2}{c|}{} & \multicolumn{3}{c|}{SNN-EP} & \multicolumn{3}{c}{DRM-M} \\ \cline{3-8}
			$k$ & ($n_1,n_2$) & \textbf{$err_{\lambda}$} & \textbf{$err_{L^2}$} & \textbf{$err_{H^1}$} & \textbf{$err_{\lambda}$} & \textbf{$err_{L^2}$} & \textbf{$err_{H^1}$} \\ \hline
			1  & (1,1) & 3.960e-15 & 2.799e-08 & 3.912e-07 & 3.868e-07 & 1.222e-04 & 7.086e-04 \\
			2  & (2,1) & 4.996e-14 & 4.711e-08 & 4.239e-07 & 1.325e-07 & 1.135e-04 & 4.374e-04 \\
			3  & (1,2) & 1.526e-13 & 5.802e-08 & 5.366e-07 & 2.186e-07 & 1.249e-04 & 6.110e-04 \\
			4  & (2,2) & 1.449e-13 & 6.954e-08 & 5.151e-07 & 7.583e-08 & 9.439e-05 & 3.828e-04 \\
			5  & (3,1) & 1.523e-13 & 8.539e-08 & 5.753e-07 & 1.126e-07 & 2.027e-04 & 5.936e-04 \\
			6  & (1,3) & 1.892e-13 & 9.419e-08 & 6.128e-07 & 1.608e-07 & 1.709e-04 & 4.993e-04 \\
			7  & (3,2) & 1.256e-13 & 1.137e-07 & 6.993e-07 & 8.341e-08 & 1.899e-04 & 4.586e-04 \\
			8  & (2,3) & 3.792e-13 & 1.284e-07 & 7.331e-07 & 1.163e-07 & 1.825e-04 & 4.520e-04 \\
			9  & (4,1) & 9.622e-14 & 1.416e-07 & 7.757e-07 & 3.046e-07 & 1.360e-03 & 1.962e-03 \\
			10 & (1,4) & 2.531e-13 & 1.436e-07 & 7.976e-07 & 1.193e-06 & 5.145e-04 & 9.427e-04 \\
			11 & (3,3) & 1.923e-13 & 1.774e-07 & 9.253e-07 & 1.106e-07 & 2.519e-04 & 4.826e-04 \\
			12 & (4,2) & 1.518e-13 & 1.737e-07 & 8.788e-07 & 2.862e-07 & 1.372e-03 & 1.908e-03 \\
			13 & (2,4) & 5.876e-13 & 2.258e-07 & 1.181e-06 & 1.039e-06 & 5.832e-04 & 9.282e-04 \\
			14 & (4,3) & 2.919e-13 & 2.558e-07 & 1.096e-06 & 3.565e-07 & 1.184e-03 & 1.516e-03 \\
			15 & (3,4) & 4.750e-13 & 3.231e-07 & 1.435e-06 & 6.469e-07 & 8.103e-04 & 1.084e-03 \\ \hline
	\end{tabular}}
\end{table}
\begin{table}[htbp]
	\centering
	\caption{The epochs needed for using SNN-EP with  different $M$ to solve the 2D Laplace eigenvalue problem.}
	\label{tab:laplace_table_epoch}
	\scalebox{0.75}{
		\begin{tabular}{lllllllllllll}
			\hline
			M & 50& 100& 150& 200&250&300&350&  400  &   450  &  500   &  550   &  600   \\ \hline
			Epoch     & 177&195&   190  &   342  &   160  & 165  &   171  &  160  &   271  &  168   &  200   &  191   \\ \hline
	\end{tabular}}
\end{table}

\subsection{Harmonic oscillator problems}
In this subsection, we study the following 2D harmonic oscillator: Find $(\lambda,u)\in \mathbb{R} \times H_0^1(\rm\Omega)$ such that $\int_{\rm\Omega}\vert u\vert^2 d\rm\Omega=1$  and
$$
-\frac{1}{2}{\rm\Delta} u +Vu= \lambda u \quad \text{in} \quad \rm\Omega,
$$
where ${\rm\Omega}=\mathbb{R}^2,V(x_1,x_2)=\frac{1}{2}(a_{11}x_1^2+2a_{12}x_1x_2+x_{22}^2)$.

Two concrete cases are considered here. For the first case, we use the decoupled harmonic oscillator potential which is most common used, i.e., $a_{11}=a_{22}=1, a_{12}=0$. For the second case, we adopt a coupled harmonic oscillator potential which is used in \cite{PhysRevA.} and maintain the same accuracy of the coefficients $a_{11},a_{12},a_{22}$ as \cite{PhysRevA.}, i.e. $a_{11}=0.8851,a_{12}=-0.1382,a_{22}=1.1933$.

\textbf{Case 1: decoupled case.}
The exact eigenvalues are given as follows:
$$
\lambda_{n_1,n_2} = (\frac{1}{2}+n_1)+(\frac{1}{2}+n_2), n_1, n_2 = 0, 1, \cdots,
$$
and the corresponding eigenfunctions  are
$$
u_{n_1,n_2}(x_1, x_2) =\mathcal{H}_{n_1}(x_1)e^{-x_1^2/2}\mathcal{H}_{n_2}(x_2)e^{-x_2^2/2}, n_1, n_2 = 0, 1, \cdots, 
$$
where $\mathcal{H}_n$ is Hermite polynomials \cite{atakishiev1990difference}.

For using SNN-EP to solve the 2D decoupled harmonic oscillator above, since $\rm\Omega$ is $\mathbb{R}^2$, the Hermite-Gauss integration scheme is chosen to compute the unbounded integral, with $99$ integration points in each dimension. In the training part, we take $\epsilon=10^{-2}$. In the dimensionality reduction part, we take $\gamma=10^{-10}$.

The dimension of the reduced subspace and the relative errors of the smallest 15 eigenvalues obtained by using different $M$ are presented in Figure 2. It can be observed that when $M$ is small, the dimension of the reduced subspace $K$ increases significantly with the increase of $M$. At this stage, the relative errors of the smallest $15$ eigenvalues rapidly decay to below $10^{-11}$. When $M$ reaches $900$, the increase of the dimension of the reduced subspace $K$ slows down significantly and the relative errors of the smallest $15$ eigenvalues remain almost unchanged. Therefore, we only show results obtained by setting $M=900$ for both SNN-EP and DRM-M.

\begin{figure}[htbp]
	\label{fig:hs_fig}
	\includegraphics[height=7.2cm,width=11.6505cm]{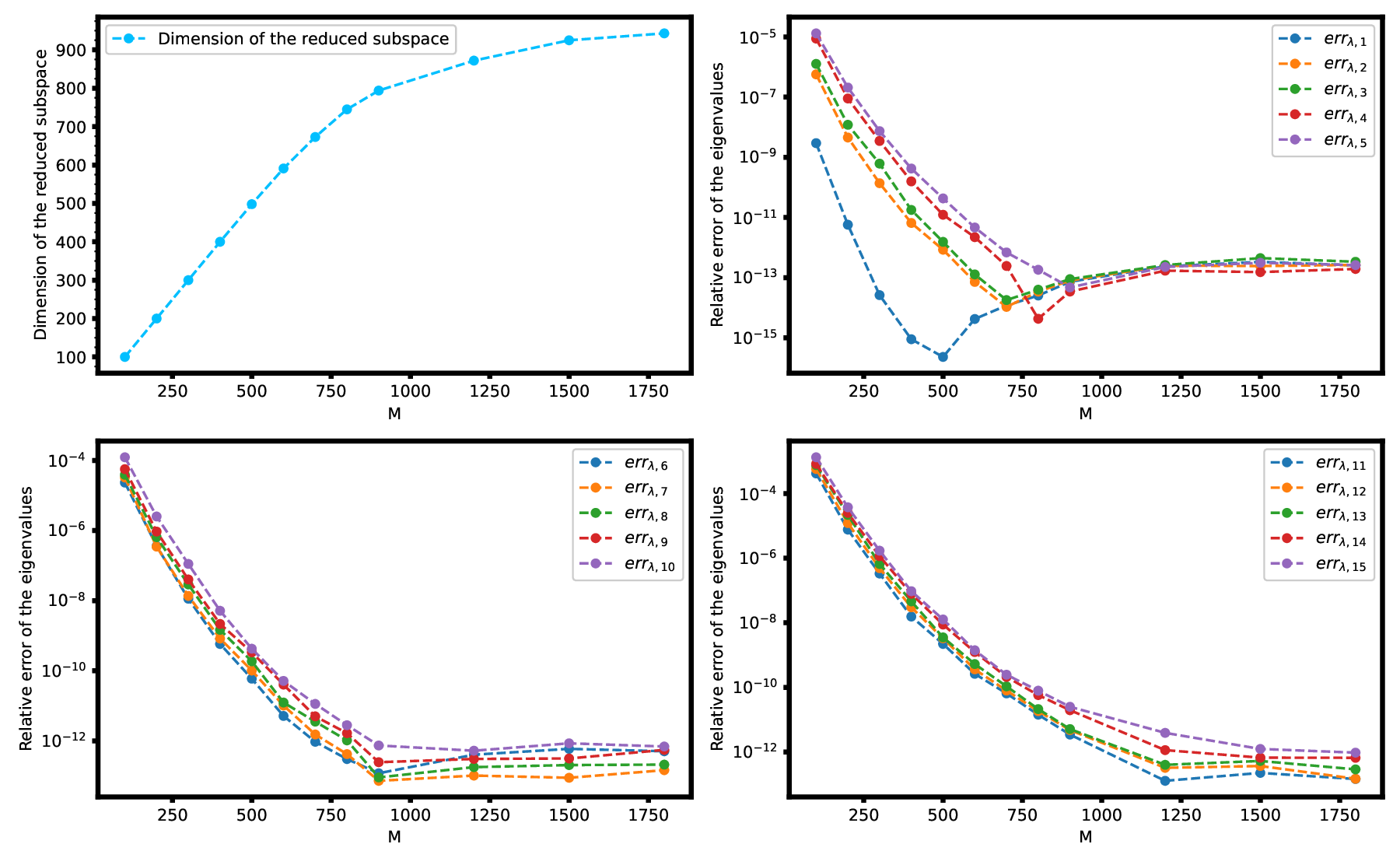}
	\centering
	\caption{The dimension of the reduced subspace and the relative errors for the 2D decoupled harmonic oscillator obtained by SNN-EP with different $M$.}
\end{figure}

The relative errors of the smallest $15$ eigenpairs obtained using SNN-EP and DRM-M with $M=900$ are presented in Table \ref{tab:hs_table}. Table \ref{tab:hs_table_epoch} provides the number of epochs needed for our method SNN-EP with various $M$. Here, it is evident that the accuracy obtained by SNN-EP is significantly superior to that obtained by DRM-M, while SNN-EP requires far fewer epochs than DRM-M. Furthermore, our results are also superior to the existing results in literature. For example, an accuracy of $10^{-2}$ is achieved using machine learning methods based on Monte Carlo integration, as demonstrated in \cite{PhysRevA.}. Additionally, our accuracy surpasses that of \cite{WANG2024112928}, which employs numerical integration based on TNN to achieve high-precision calculations. Compared to \cite{WANG2024112928}, which requires $500,000$ epochs, our method achieves the desired accuracy with only $20$ epochs.

\begin{table}[htbp]
	\centering
	\caption{The relative errors of the smallest $15$ eigenvalues of the 2D decoupled harmonic oscillator obtained by setting $M=900$.}
	\label{tab:hs_table}
	\scalebox{0.75}{
		\begin{tabular}{ll|lll|lll}
			\hline
			\multicolumn{2}{c|}{} & \multicolumn{3}{c|}{SNN-EP} & \multicolumn{3}{c}{DRM-M} \\ \cline{3-8}
			$k$ & ($n_1,n_2$) & \textbf{$err_{\lambda}$} & \textbf{$err_{L^2}$} & \textbf{$err_{H^1}$} & \textbf{$err_{\lambda}$} & \textbf{$err_{L^2}$} & \textbf{$err_{H^1}$} \\ \hline
			1  & (0,0) & 7.294e-14 & 3.977e-08 & 2.373e-07 & 2.449e-09 & 2.062e-05 & 3.913e-05 \\
			2  & (0,1) & 8.005e-14 & 6.304e-08 & 3.012e-07 & 4.052e-09 & 3.637e-05 & 6.912e-05 \\
			3  & (1,0) & 8.926e-14 & 6.095e-08 & 2.963e-07 & 8.495e-09 & 5.009e-05 & 9.015e-05 \\
			4  & (0,2) & 3.449e-14 & 1.184e-07 & 4.473e-07 & 3.324e-09 & 4.887e-05 & 7.840e-05 \\
			5  & (1,1) & 4.737e-14 & 9.585e-08 & 3.648e-07 & 4.008e-09 & 5.239e-05 & 8.106e-05 \\
			6  & (2,0) & 1.198e-13 & 9.863e-08 & 4.040e-07 & 6.147e-09 & 4.121e-05 & 6.859e-05 \\
			7  & (0,3) & 7.283e-14 & 3.044e-07 & 1.035e-06 & 3.407e-09 & 7.011e-05 & 1.035e-04 \\
			8  & (1,2) & 9.104e-14 & 1.458e-07 & 5.599e-07 & 3.609e-09 & 6.003e-05 & 8.251e-05 \\
			9  & (2,1) & 2.474e-13 & 1.632e-07 & 5.801e-07 & 7.057e-09 & 6.483e-05 & 9.069e-05 \\
			10 & (3,0) & 7.412e-13 & 2.373e-07 & 8.030e-07 & 8.040e-09 & 7.635e-05 & 1.054e-04 \\
			11 & (1,3) & 3.457e-12 & 1.312e-06 & 3.803e-06 & 6.594e-09 & 1.537e-04 & 1.734e-04 \\
			12 & (2,2) & 4.765e-12 & 9.109e-07 & 2.550e-06 & 7.201e-09 & 1.784e-04 & 2.085e-04 \\
			13 & (3,1) & 5.116e-12 & 1.221e-06 & 3.586e-06 & 1.029e-08 & 1.329e-04 & 1.632e-04 \\
			14 & (0,4) & 1.962e-11 & 1.923e-06 & 5.044e-06 & 1.426e-08 & 1.176e-04 & 1.424e-04 \\
			15 & (4,0) & 2.506e-11 & 1.185e-06 & 3.308e-06 & 2.001e-08 & 1.511e-04 & 1.772e-04 \\ \hline
	\end{tabular}}
\end{table}


\begin{table}[htbp]
	\centering
	\caption{The epochs needed for using SNN-EP with different $M$ to solve the 2D decoupled harmonic oscillator.}
	\label{tab:hs_table_epoch}
	\scalebox{0.75}{
		\begin{tabular}{lllllllllllll}
			\hline
			M & 100 & 200 & 300 & 400 & 500 & 600 &   700  &  800  &   900  &  1200   &  1500   &  1800   \\ \hline
			Epoch     &   20  &  21   &   23  &   23  &   22  & 20  &   20  &  21  &   21  &  21   &  20   &  22   \\ \hline
	\end{tabular}}
\end{table}

\textbf{Case 2: coupled case.}
The coupled harmonic oscillator defined by the above potential energy term can be decoupled into the following diagonal form:
$$V(y_1,y_2)=\frac{1}{2}(\mu_1y_1^2+\mu_2y_2^2),$$
where $y_1= -0.9339352418x_1 - 0.3574422527x_2$ and $y_2= 0.3574422527x_1 - 0.93393524$ $18x_2$, $\mu_1=0.8322071257$, $\mu_2=1.2461928742$.

Based on the decoupled form, the exact eigenvalues are given as follows:
$$
\lambda_{n_1,n_2} = (\frac{1}{2}+n_1)\mu_1^{1/2}+(\frac{1}{2}+n_2)\mu_2^{1/2}, n_1, n_2 = 0, 1, \cdots,
$$
and the corresponding eigenfunctions  are
$$
u_{n_1,n_2}(x_1, x_2) =\mathcal{H}_{n_1}(\mu_1^{1/4}y_1)e^{-\mu_1^{1/2}y_1^2/2}\mathcal{H}_{n_2}(\mu_2^{1/4}y_2)e^{-\mu_2^{1/2}y_2^2/2}, n_1, n_2 = 0, 1, \cdots.
$$

For using SNN-EP to solve the 2D coupled harmonic oscillator above, since $\rm\Omega$ is $\mathbb{R}^2$, we choose the Hermite-Gauss integration scheme, with $99$ integration points in each dimension. In the training part, we take $\epsilon=10^{-2}$. In the dimensionality reduction part, we take $\gamma=10^{-10}$.

The dimension of the reduced subspace and the relative errors of the smallest $15$ eigenvalues obtained by SNN-EP with different $M$ are presented in Figure 3. A numerical performance similar to the decoupled harmonic oscillator can be observed. Accordingly, we present detailed results for $M=900$ for both SNN-EP and DRM-M.

\begin{figure}[htbp]
	\label{fig:chs_fig}
	\includegraphics[height=7.2cm,width=11.6505cm]{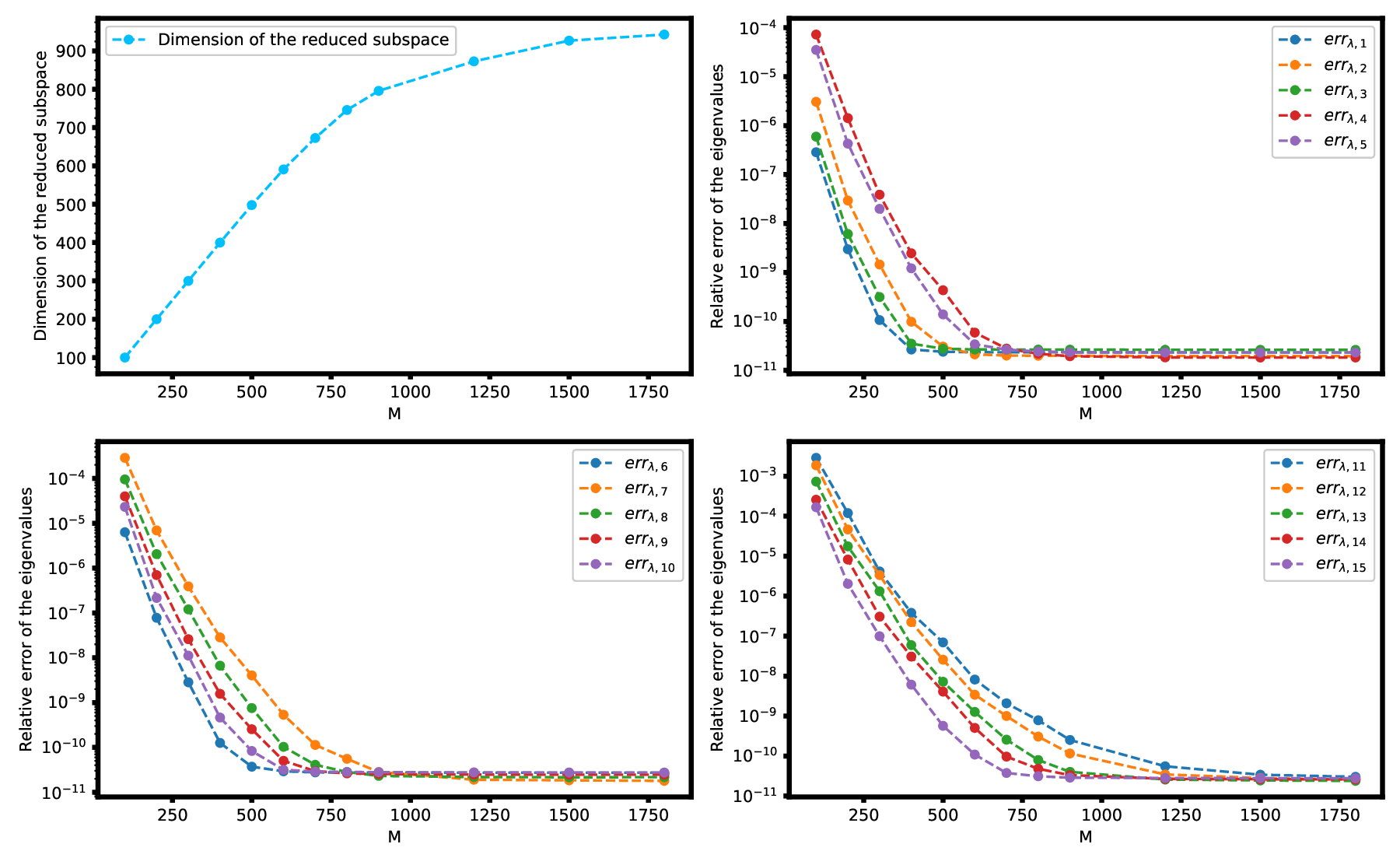}
	\centering
	\caption{The dimension of the reduced subspace and the relative errors of eigenvalues of the 2D coupled harmonic oscillator obtained by SNN-EP with different $M$.}
\end{figure}

The relative errors of the smallest $15$ eigenpairs of  the 2D coupled harmonic oscillator obtained by using  SNN-EP and DRM-M  with $M=900$ are shown in Table \ref{tab:chs_table}. Similar to the decoupled harmonic oscillator, our method obtains approximate solutions with higher accuracy while requiring fewer epochs during the training, as shown in Table \ref{tab:chs_table_epoch}, compared to  DRM-M and other existing results in \cite{PhysRevA.} and \cite{WANG2024112928}.

\begin{table}[htbp]
	\centering
	\caption{The relative errors of the smallest $15$ eigenvalues of the 2D coupled harmonic oscillator obtained by setting $M=900$.}
	\label{tab:chs_table}
	\scalebox{0.75}{
		\begin{tabular}{ll|lll|lll}
			\hline
			\multicolumn{2}{c|}{} & \multicolumn{3}{c|}{SNN-EP} & \multicolumn{3}{c}{DRM-M} \\ \cline{3-8}
			$k$ & ($n_1,n_2$) & \textbf{$err_{\lambda}$} & \textbf{$err_{L^2}$} & \textbf{$err_{H^1}$} & \textbf{$err_{\lambda}$} & \textbf{$err_{L^2}$} & \textbf{$err_{H^1}$} \\ \hline
			1  & (0,0) & 2.330e-11 & 4.430e-08 & 2.520e-07 & 1.154e-07 & 1.361e-04 & 2.723e-04 \\
			2  & (1,0) & 1.979e-11 & 6.415e-08 & 3.183e-07 & 2.964e-08 & 9.442e-05 & 1.736e-04 \\
			3  & (0,1) & 2.647e-11 & 7.323e-08 & 3.517e-07 & 2.188e-08 & 9.079e-05 & 1.600e-04 \\
			4  & (2,0) & 1.939e-11 & 2.710e-07 & 1.011e-06 & 6.870e-08 & 1.912e-04 & 2.966e-04 \\
			5  & (1,1) & 2.348e-11 & 1.704e-07 & 6.191e-07 & 1.830e-08 & 9.054e-05 & 1.414e-04 \\
			6  & (0,2) & 2.740e-11 & 1.188e-07 & 4.584e-07 & 1.826e-08 & 1.172e-04 & 1.754e-04 \\
			7  & (3,0) & 2.829e-11 & 1.107e-06 & 3.355e-06 & 5.469e-08 & 1.966e-04 & 2.674e-04 \\
			8  & (2,1) & 2.312e-11 & 4.443e-07 & 1.350e-06 & 4.192e-08 & 2.110e-04 & 2.704e-04 \\
			9  & (1,2) & 2.545e-11 & 2.514e-07 & 8.112e-07 & 1.655e-08 & 1.173e-04 & 1.616e-04 \\
			10 & (0,3) & 2.798e-11 & 1.478e-07 & 5.192e-07 & 1.759e-08 & 1.341e-04 & 1.743e-04 \\
			11 & (4,0) & 2.524e-10 & 6.010e-06 & 1.645e-05 & 9.014e-08 & 3.495e-04 & 4.336e-04 \\
			12 & (3,1) & 1.159e-10 & 4.032e-06 & 1.036e-05 & 6.991e-08 & 3.703e-04 & 4.364e-04 \\
			13 & (2,2) & 4.021e-11 & 1.626e-06 & 4.390e-06 & 3.048e-08 & 2.252e-04 & 2.719e-04 \\
			14 & (1,3) & 3.324e-11 & 1.120e-06 & 2.951e-06 & 4.191e-08 & 2.746e-04 & 3.394e-04 \\
			15 & (0,4) & 2.874e-11 & 3.594e-07 & 1.063e-06 & 3.972e-08 & 2.894e-04 & 3.600e-04 \\ \hline
	\end{tabular}}
\end{table}

\begin{table}[htbp]
	\centering
	\caption{The epochs needed for using SNN-EP with different $M$ to deal with the  2D coupled harmonic oscillator.}
	\label{tab:chs_table_epoch}
	\scalebox{0.75}{
		\begin{tabular}{lllllllllllll}
			\hline
			M & 100 & 200 & 300 & 400 & 500 & 600 &   700  &  800  &   900  &  1200   &  1500   &  1800   \\ \hline
			Epoch     &   20  &  21   &   22  &   23  &   22  & 20  &   20  &  21  &   20  &  21   &  20   &  23   \\ \hline
	\end{tabular}}
\end{table}

\subsection{Hydrogen atom}
We consider the following Schrodinger equation for  Hydrogen atom in 2D: Find $(\lambda,u)\in \mathbb{R} \times H_0^1(\rm\Omega)$ such that $\int_{\rm\Omega}\vert u\vert^2 d\rm\Omega=1$  and
$$
-\frac{1}{2}{\rm\Delta} u +Vu= \lambda u \quad \text{in} \quad \rm\Omega,
$$
where ${\rm\Omega}=\mathbb{R}^2,V(x_1,x_2)=-\frac{1}{\sqrt{x_1^2+x_2^2}}$.

The exact eigenvalues are as follows:
$$
\lambda_{n_1,n_2} = -\frac{2}{(2|n_1|+2n_2+1)^2}, |n_1| , n_2 = 0, 1, \cdots,,
$$
and the corresponding eigenfunctions in spherical coordinates are
$$
u_{n_1,n_2}(r, \theta) =e^{in_1\theta}e^{-\frac{2}{2|n_1|+2n_2+1} r}r^{|n_1|}\mathcal{P}_{|n_1|,n_2}(r), |n_1| , n_2 = 0, 1, \cdots, 
$$
where $\mathcal{P}_{m,n}(r)=\sum\limits_{i=0}^n a_i r^i,$ $a_{i+1}=\frac{\beta(2m+2i+1)-2}{(i+1)(i+2m+1)}a_i$, $\beta=\frac{2}{2|m|+2n+1} ,i=0,1,\cdots,n-1$ and $a_0$ is chosen to normalize $u_{n_1,n_2}(r, \theta)$.

To deal with the singularity of the potential function $V$, we generate numerical integration points in spherical coordinates, i.e.,  $\int_{\mathbb{R}^2}g(x_1,x_2)dx_1dx_2=\int_{0}^{2\pi}\int_{0}^{+\infty}g(r\cos\theta,$ $r\sin\theta)r^2sin\theta drd\theta,$ where $g$ is a general function.

For the integration along  $r$, we choose the Hermite-Gauss integration scheme, with 99 integration points. {For the integration along  $\theta$, we choose the Legendre-Gauss integration scheme, with 4 integration points in each subinterval and  the number of subintervals $20$ . In training part, we take $\epsilon=10^{-2}$. In the dimensionality reduction part, we take $\gamma=10^{-10}$.

	Different from the previous examples, here, we apply 5 different modification
	functions $f_0,\cdots,f_4$ to each neural network-based function $\varphi_j$. Then we train the parameters
	and perform dimensionality reduction operation to the modified functions $\{f_i\varphi_j\}_{i=0,j=1}^{4,M}$, here $f_i(x_1,x_2) = e^{-\frac{2}{2i+1}{\sqrt{x_1^2+x_2^2}}},i = 0,1,2,3.$

	The dimension of the reduced subspace and the relative errors of the smallest 15 eigenvalues obtained by using different $M$ are presented in Figure 4. It can be observed that when $M$ is small, the dimension of the reduced subspace $K$ increases significantly with the increase of $M$. At this stage, the relative errors of the smallest $15$ eigenvalues rapidly decay to below $10^{-12}$. When $M$ reaches $600$, the increase of the dimension of the reduced subspace $K$ slows down significantly. When $M$ reaches $1000$, the relative errors of the smallest $15$ eigenvalues remain almost unchanged. Therefore, we only show results obtained by setting $M=1000$ for both SNN-EP and DRM-M in Table \ref{tab:hy_table}. From these results, we can see clearly that the accuracy obtained by SNN-EP is significantly superior to that obtained by DRM-M, while SNN-EP requires much fewer epochs as listed in Table \ref{tab:hy_table_epoch}.
	
	It should be pointed out that the selection of the above $f_i$ is too specific. In fact, we can also choose $f_i$ in the following more general way to obtain $\{f_i\varphi_j\}_{i=1,j=1}^{N,M}$, e.g. $f_i(\bm x)=e^{-\beta_i{\Vert x\Vert}}$, $i=1,\cdots,N.$ Here, $\beta_i$ is randomly chosen from a uniform distribution and $N$ represents the number of randomly selected $\beta_i$. For the randomly selected $\beta_i$, we have also done some numerical experiments, whose results are presented in Figure 5. These results show that as long as we use enough number of $\beta_i$, even random ones can achieve a similar effect, and for this problem, $N$=10 is sufficient.
	
	\begin{figure}[htbp]
			\label{fig:hy_fig}
		\includegraphics[height=7.2cm,width=11.6505cm]{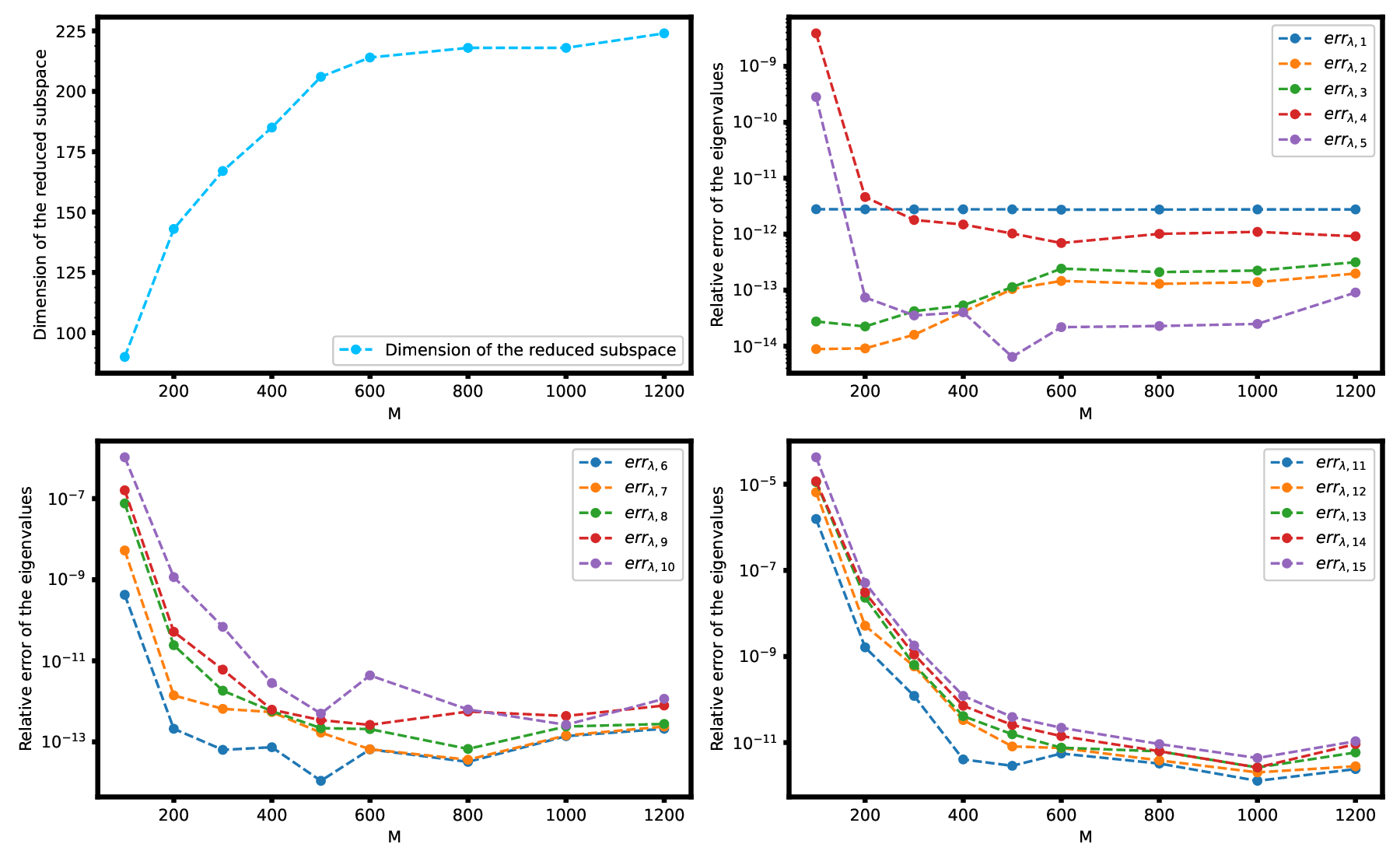}
		\centering
		\caption{The dimension of the reduced subspace and the relative error of eigenvalues of the Schrodinger equation for  Hydrogen atom obtained by SNN-EP with different $M$.}
	\end{figure}
	
	\begin{table}[htbp]
		\centering
		\caption{The relative errors of the smallest $15$ eigenvalues of the Schrodinger equation for  Hydrogen atom  obtained by SNN-EP with  $M=1000$.}
		\label{tab:hy_table}
		\scalebox{0.75}{
			\begin{tabular}{ll|lll|lll}
				\hline
				\multicolumn{2}{c|}{} & \multicolumn{3}{c|}{SNN-EP} & \multicolumn{3}{c}{DRM-M} \\ \cline{3-8} 
				$k$ & ($n_1,n_2$) & \textbf{$err_{\lambda}$} & \textbf{$err_{L^2}$} & \textbf{$err_{H^1}$} & \textbf{$err_{\lambda}$} & \textbf{$err_{L^2}$} & \textbf{$err_{H^1}$} \\ \hline
				1  & ( 0,0) & 2.770e-12 & 6.092e-09 & 1.642e-08 & 2.847e-07 & 1.144e-04 & 1.370e-04 \\
				2  & ( 0,1) & 1.392e-13 & 5.620e-08 & 3.046e-07 & 6.322e-04 & 1.047e-02 & 5.763e-02 \\
				3  & ( 1,0) & 2.244e-13 & 6.228e-08 & 8.909e-08 & 1.479e-05 & 2.991e-03 & 6.391e-03 \\
				4  & (-1,0) & 1.101e-12 & 4.646e-08 & 5.109e-08 & 6.761e-05 & 5.415e-03 & 1.224e-02 \\
				5  & ( 0,2) & 2.498e-14 & 6.278e-08 & 1.481e-07 & 5.019e-05 & 5.351e-03 & 3.057e-02 \\
				6  & ( 1,1) & 1.381e-13 & 1.641e-07 & 6.657e-07 & 5.285e-05 & 6.628e-03 & 1.612e-02 \\
				7  & (-1,1) & 1.313e-13 & 1.148e-07 & 3.718e-06 & 7.372e-05 & 7.391e-03 & 1.764e-02 \\
				8  & ( 2,0) & 2.392e-13 & 7.195e-08 & 1.482e-07 & 1.530e-04 & 8.331e-03 & 3.442e-02 \\
				9  & (-2,0) & 4.348e-13 & 5.816e-08 & 1.294e-07 & 2.125e-04 & 4.244e-03 & 1.203e-02 \\
				10 & ( 0,3) & 2.633e-13 & 4.133e-07 & 1.408e-06 & 6.733e-04 & 3.532e-02 & 1.011e-01 \\
				11 & ( 1,2) & 1.293e-12 & 4.411e-07 & 7.443e-07 & 7.967e-04 & 5.432e-02 & 1.230e-01 \\
				12 & (-1,2) & 2.032e-12 & 2.933e-07 & 5.402e-07 & 1.243e-03 & 8.099e-02 & 2.104e-01 \\
				13 & ( 2,1) & 2.639e-12 & 5.727e-07 & 1.371e-06 & 1.452e-03 & 4.912e-02 & 2.078e-01 \\
				14 & (-2,1) & 2.640e-12 & 6.261e-07 & 1.510e-06 & 1.524e-03 & 2.474e-02 & 6.840e-02 \\
				15 & ( 3,0) & 4.376e-12 & 4.827e-07 & 8.022e-07 & 2.384e-03 & 2.179e-02 & 7.715e-02 \\ \hline
		\end{tabular}}
	\end{table}
	
	\begin{table}[htbp]
		\centering
		\caption{The epochs needed for using SNN-EP with different $M$ to solve the Schrodinger equation for  Hydrogen atom.}
		\label{tab:hy_table_epoch}
		\scalebox{0.75}{
			\begin{tabular}{llllllllll}
				\hline
				M & 100 & 200 & 300 & 400 & 500 & 600   &  800  &   1000  &  1200         \\ \hline
				Epoch     &   129  &  91   &   167  &   46  &   127  & 81  &   178  &  119  &   246    \\ \hline
		\end{tabular}}
	\end{table}

	\begin{figure}[htbp]
		\label{fig:hy_random_fig}
		\includegraphics[height=7.2cm,width=11.6505cm]{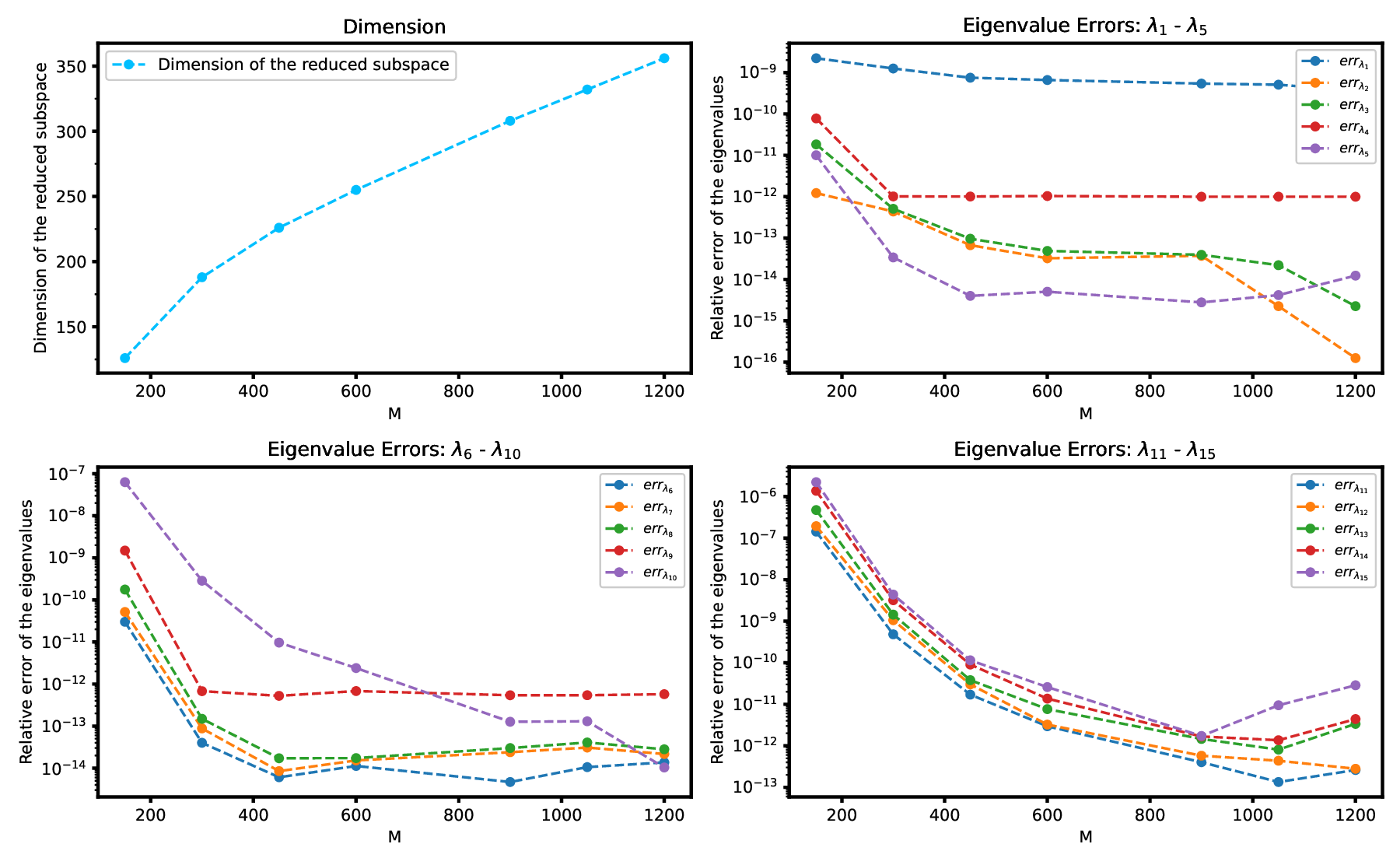}
		\centering
		\caption{The dimension of the reduced subspace and the relative error of eigenvalues of the Schrodinger equation for  Hydrogen atom obtained by SNN-EP (10 random $\beta_i$) with different $M$.}
	\end{figure}

	\section{Conclusions}
	\label{sec:conclusions}
	In this paper, we have proposed a subspace method based on neural networks for eigenvalue problems. This method can in fact be viewed as a method combining the machine learning with the classical methods. That is, first use the machine learning method and model order reduction method to get some neural network-based orthogonal basis, then calculate the Galerkin projection of the eigenvalue problem onto the subspace spanned by the orthogonal basis and obtain an approximate solution, just like what the classical methods do. Our method can achieve high precision with only a few training epochs. In fact, from the numerical experiments, we can see that our method can obtain approximate eigenpairs with the error of eigenvalues being below $10^{-11}$ and the $L^2$ norm and $H^1$ norm  error of eigenfunctions being nearly below $10^{-7}$, but with only tens or hundreds of epochs.

\section*{Appendix A. The impact of training on solution accuracy}
In this section, we provide some numerical examples to illustrate how the parameters $\theta$ used in the training affect the solution accuracy. After each gradient descent step, we solve the eigenvalue problem in the subspace spanned by the basis functions $\psi_1(\bm x,\theta),\cdots,\psi_K(\bm x,\theta)$ which are obtained by the newly updated parameters $\theta$ and model order reduction. We take the results of 2D Laplace eigenvalue problem on a square domain $[0,1]^2$ as an example, whose details are introduced in Section 4. The number of neurons of the subspace layer $M$ is set to 300, and all other parameters are the same with Section 4.  Figure 6 shows the solution accuracy at different gradient descent steps. We can see that for the smallest few eigenvalues, it is easy to obtain high accuracy approximation without training, while for the next smallest eigenvalues, training can significantly improve the accuracy of the results.
\begin{figure}[htbp]
	\label{fig:train_fig}
	\includegraphics[height=7.2cm,width=11.6505cm]{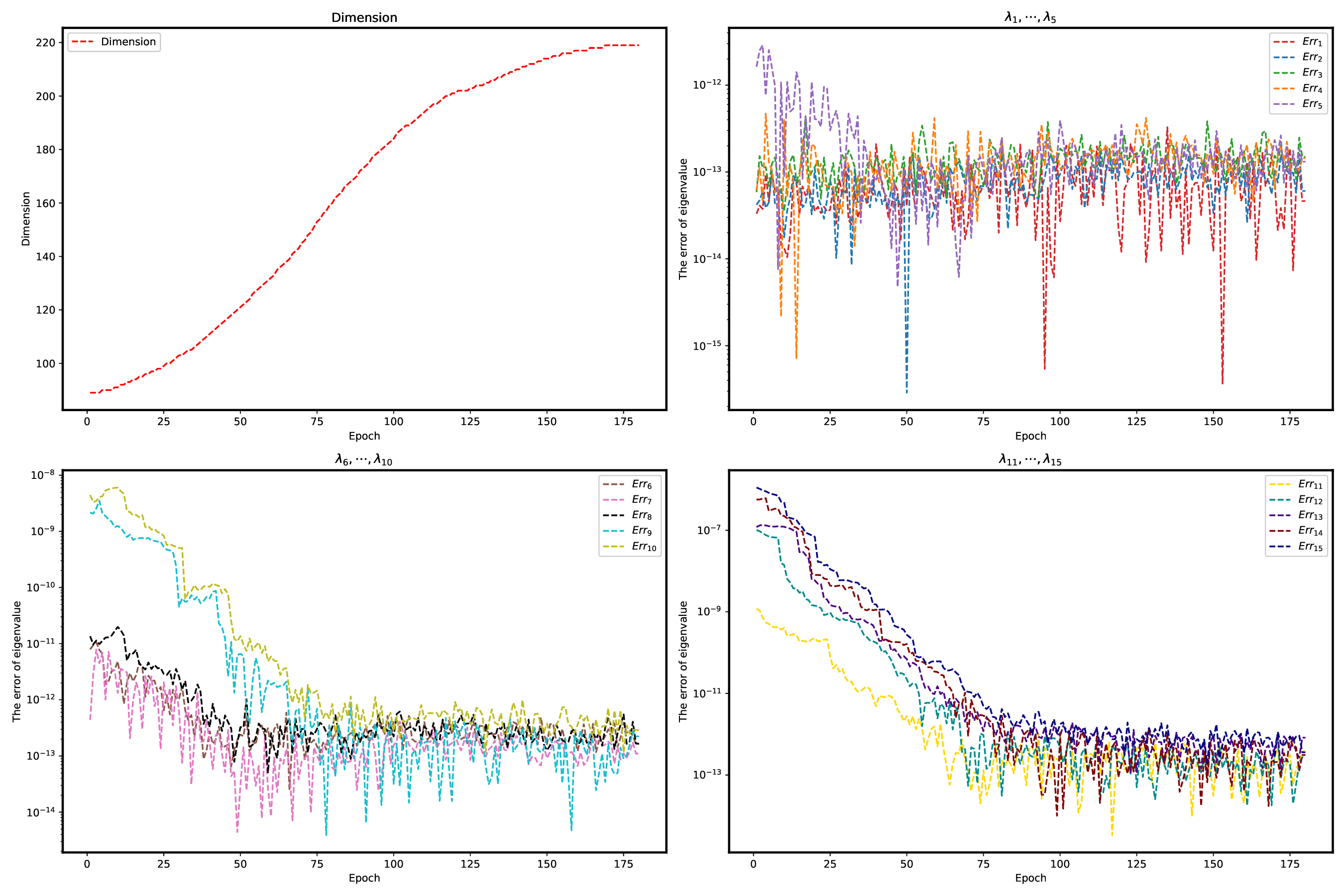}
	\centering
	\caption{The effect of trained parameters $\theta$ on the accuracy of solution.}
\end{figure}

Furthermore, we see the dimension of the subspace spanned by the basis $\psi_1(\bm x,\theta),$ $\cdots,$ $\psi_K(\bm x,\theta)$ generated by training and model order reduction continues to increase as the training steps increase. Consequently, more information is preserved during the dimensionality reduction part, which helps to explain the increase in solution accuracy throughout the training process.

\section*{Appendix B. The necessity of employing dimensionality reduction technique}
In this section, we provide some numerical tests that can show the necessity of employment dimensionality reduction technique. In fact, we have applied the version of Algorithm 3 that excludes the dimensionality reduction step to solution of 2D Laplace eigenvalue problem on a square domain $[0,1]^2$ as introduced in Section 4, aiming to solve the first 15 eigenpairs. In this situation, the algebraic eigenvalue problem (\ref{MEP1}) is extremely ill-conditioned. Here, we adopt the solver $scipy.linalg.eig$ in Python, which can deal with the ill-conditioned systems  in some sense. All other settings are exactly the same as in Section 4. Numerical experiments for different $M$ have been carried out and the results are presented in Table \ref{tab:2DL_N}.

\begin{table}[htbp]
	\centering
	\caption{The relative errors of the smallest $15$ eigenvalues obtained by the version of Algorithm 3 that excludes the dimensionality reduction step under different $M$ for the 2D Laplace eigenvalue problem on a square domain $[0,1]^2$.}
	\label{tab:2DL_N}
	\scalebox{0.75}{
		\begin{tabular}{lllllll}
			\hline
			Dimension & 100 & 200 & 300 & 400 & 500&600  \\ \hline
			$err_{\lambda,1}$&3.901e-12&1.315e-11&3.477e-04&2.236e-04&4.817e-10&3.047e-04\\
			$err_{\lambda,2}$&1.586e-11&2.618e-09&9.126e-05&5.287e-05&1.593e-08&3.785e-10\\
			$err_{\lambda,3}$&5.342e-11&4.210e-03&1.255e-03&1.851e-02&1.476e-01&5.980e-02\\
			$err_{\lambda,4}$&1.197e-09&6.017e-04&1.471e-06&3.563e-04&1.745e-04&5.853e-04\\
			$err_{\lambda,5}$&9.387e-10&5.089e-12&8.218e-06&4.899e-03&1.009e-03&8.677e-03\\
			$err_{\lambda,6}$&5.093e-09&7.511e-11&9.078e-03&1.221e-02&1.116e-01&1.162e-02\\
			$err_{\lambda,7}$&5.447e-09&3.872e-10&8.271e-06&7.278e-02&3.229e-04&3.899e-03\\
			$err_{\lambda,8}$&8.941e-09&3.102e-01&6.478e-03&2.485e-01&2.230e-03&4.813e-03\\
			$err_{\lambda,9}$&3.880e-08&5.917e-02&3.574e-04&9.198e-03&4.877e-03&6.637e-03\\
			$err_{\lambda,10}$&8.691e-08&1.811e-01&4.090e-02&1.514e-02&2.628e-02&2.852e-02\\
			$err_{\lambda,11}$&6.166e-08&3.851e-01&2.117e-03&4.508e-03&2.238e-02&9.721e-02\\
			$err_{\lambda,12}$&2.687e-07&3.014e-01&1.180e-03&1.608e-03&4.303e-02&1.707e-02\\
			$err_{\lambda,13}$&4.700e-07&4.493e-01&5.860e-02&9.130e-03&5.374e-02&4.256e-02\\
			$err_{\lambda,14}$&3.048e-07&2.871e-01&5.064e-03&9.433e-03&6.394e-03&3.686e-03\\
			$err_{\lambda,15}$&7.798e-07&4.797e-01&8.226e-02&9.797e-03&4.927e-02&1.188e-02\\ \hline
	\end{tabular}}
\end{table}

From Table \ref{tab:2DL_N} we can see that when $M$ is small, our method can   obtain approximate eigenvalues with high accuracy. However, with the increase of $M$, the accuracy of eigenvalues is not improved  but rapidly deteriorates, contrary to our expectations.

After some tests, we found that the main reason for the poor accuracy may lie in the extreme ill-conditioning of the matrices $A$ and $B$ in (\ref{MEP1}). From Table \ref{tab:2DL_N1}, where the condition numbers of the matrices $A$ and $B$ corresponding to different $M$ are shown, we can see that as the increase of $M$, the condition numbers of matrices $A$ and $B$ increase a lot. By taking a detailed look at  Table \ref{tab:2DL_N} and Table \ref{tab:2DL_N1}, we believe that it is the increasingly ill-conditioned matrices $A$ and $B$  that prevent the solver from obtaining solutions with high accuracy. In fact, a similar phenomenon has also been observed in random feature method \cite{chen2024optimization}, and the low-rank structure of  random feature matrix has been systematically studied in \cite{alaoui2015fast,bach2013sharp,rudi2015less,williams2000using}. 

\begin{table}[htbp]
	\centering
	\renewcommand{\arraystretch}{1.2}
	\caption{The condition number $\kappa(A)$ and $\kappa(B)$ generated by the version of Algorithm 3 that excludes the dimensioality reduction step under different $M$ for the 2D Laplace eigenvalue problem on a square domain $[0,1]^2$.}
	\label{tab:2DL_N1}
	\scalebox{0.75}{
		\begin{tabular}{lllllll}
			\hline
			M & 100 & 200 & 300 & 400 & 500 &600 \\ \hline
			$\kappa(A)$&7.945e+15 & 2.383e+21 & 1.754e+26 & 8.161e+30 & 3.256e+33 & 6.145e+33\\
			$\kappa( B)$&2.954e+17 & 4.261e+23 & 2.024e+28 & 6.245e+32 & 8.485e+33 & 2.101e+34\\\hline
	\end{tabular}}
\end{table}
Hence, these results show the necessity of employing some dimensionality reduction technique in Algorithm \ref{alg:A3} to construct neural network-based orthogonal basis functions.

In fact, we have done a series of tests, which show that this step can improve the condition numbers of the stiffness matrix and the mass matrix greatly. The condition numbers of the stiffness matrix and the mass matrix with the dimensionality reduction operation with $\gamma=10^{-9},10^{-10}$ and $10^{-11}$ are presented in Table \ref{tab:2DL_N2}. Comparing Table \ref{tab:2DL_N2} with Table \ref{tab:2DL_N1}, we see clearly that the condition number of the matrices obtained without the dimensionality reduction operation increases rapidly with the increase of $M$. In contrast, the condition numbers of the matrices obtained with the dimensionality reduction operation can be controlled within an acceptable range and do not increase with $M$. 

\begin{table}[htbp]
\centering
\caption{The condition number $\kappa(A)$ and $\kappa(B)$ obtained by Algorithm \ref{alg:A3} with different $M$ and dimensionality reduction parameter $\gamma$ for the 2D Laplace eigenvalue problem on a square domain $[0,1]^2$.}
\label{tab:2DL_N2}
\scalebox{0.75}{
	\renewcommand{\arraystretch}{1.2}
	\begin{tabular}{l|lll|lll}
		\hline
		\multicolumn{1}{c|}{$M$} & \multicolumn{3}{c|}{\textbf{$\kappa(A)$}} & \multicolumn{3}{c}{\textbf{$\kappa(B)$}} \\ \cline{2-7} 
		& $\gamma=10^{-11}$ & $\gamma=10^{-10}$ & $\gamma=10^{-9}$ & $\gamma=10^{-11}$ & $\gamma=10^{-10}$ & $\gamma=10^{-9}$ \\ \hline
		100 & 206.26 & 206.26 & 199.82 & 1.00 & 1.00 & 1.00 \\
		200 & 513.95 & 463.29 & 353.87 & 1.00 & 1.00 & 1.00 \\
		300 & 666.72 & 540.08 & 369.34 & 1.00 & 1.00 & 1.00 \\
		400 & 742.69 & 542.86 & 376.37 & 1.00 & 1.00 & 1.00 \\
		500 & 659.94 & 468.31 & 374.71 & 1.00 & 1.00 & 1.00 \\
		600 & 819.00 & 589.94 & 451.74 & 1.00 & 1.00 & 1.00 \\ \hline
\end{tabular}}
\end{table}

\bibliographystyle{siam}
\bibliography{references}

\end{document}